# MOCCA: mirrored convex/concave optimization for nonconvex composite functions

Rina Foygel Barber and Emil Y. Sidky

04.23.16

**Abstract**

Many optimization problems arising in high-dimensional statistics decompose naturally into a sum of several terms, where the individual terms are relatively simple but the composite objective function can only be optimized with iterative algorithms. In this paper, we are interested in optimization problems of the form $\mathsf{F}(Kx) + \mathsf{G}(x)$, where $K$ is a fixed linear transformation, while $\mathsf{F}$ and $\mathsf{G}$ are functions that may be nonconvex and/or nondifferentiable. In particular, if either of the terms are nonconvex, existing alternating minimization techniques may fail to converge; other types of existing approaches may instead be unable to handle nondifferentiability. We propose the MOCCA (mirrored convex/concave) algorithm, a primal/dual optimization approach that takes local convex approximation to each term at every iteration. Inspired by optimization problems arising in computed tomography (CT) imaging, this algorithm can handle a range of nonconvex composite optimization problems, and offers theoretical guarantees for convergence when the overall problem is approximately convex (that is, any concavity in one term is balanced out by convexity in the other term). Empirical results show fast convergence for several structured signal recovery problems.

## 1 Introduction

We consider the problem of minimizing a composite objective function of the form

$$\mathsf{F}(Kx) + \mathsf{G}(x) \tag{1}$$

over $x \in \mathbb{R}^d$, where $K \in \mathbb{R}^{m\times d}$ is a fixed linear operator, and $\mathsf{F}$ and $\mathsf{G}$ are functions which are potentially nonconvex and/or nondifferentiable. Optimization problems of this form arise in many applications, and in particular, the algorithm developed here was motivated by an image reconstruction problem for computed tomography (CT), an imaging technology used often in medicine and in other domains.

When $\mathsf{F}$ and $\mathsf{G}$ are both convex, many existing methods are well-equipped to handle this optimization problem, even in high dimensions. For example, the alternating direction method of multipliers (ADMM) [4] and related primal/dual methods yield effective algorithms when the functions $\mathsf{F}$ and $\mathsf{G}$ both have inexpensive proximal maps, defined as

$$\mathsf{Prox}_{\mathsf{F}}(u) = \arg\min_{w}\left\{\frac{1}{2}\|w-u\|^2 + \mathsf{F}(w)\right\}$$

for any $u \in \mathbb{R}^m$, and same for $\mathsf{Prox}_{\mathsf{G}}$ defined on $\mathbb{R}^d$. Primal/dual methods, such as the alternating direction method of multipliers (ADMM) [4], are especially effective if $\mathsf{F}$ has an inexpensive proximal map but the linear transformation of the same function, i.e. the function $x \mapsto \mathsf{F}(Kx)$, does not. If $\mathsf{F}$ does not offer an inexpensive proximal map but is instead smoothly differentiable (even if nonconvex), while $\mathsf{G}$ does have a fast proximal map, then methods such as proximal gradient descent (see e.g. [27, 2]) can be applied instead of a primal/dual method.

In this paper, we consider a more challenging setting where $\mathsf{F}$ and $\mathsf{G}$ may both be nonconvex and nondifferentiable. For instance, we can consider working with functions that can be decomposed as $\mathsf{F} = \mathsf{F}_{\mathsf{cvx}} + \mathsf{F}_{\mathsf{diff}}$ and $\mathsf{G} = \mathsf{G}_{\mathsf{cvx}} + \mathsf{G}_{\mathsf{diff}}$, where we assume that $\mathsf{F}_{\mathsf{cvx}}, \mathsf{G}_{\mathsf{cvx}}$ are convex but do not need to be differentiable, while $\mathsf{F}_{\mathsf{diff}}, \mathsf{G}_{\mathsf{diff}}$ are potentially nonconvex but are differentiable. (If $\mathsf{F}_{\mathsf{diff}}$ and $\mathsf{G}_{\mathsf{diff}}$ are concave, then this type of optimization problem is often referred to as "convex/concave".) As we will see, this formulation arises naturally in a range of applications, but in general, cannot be handled by existing methods that are designed with convexity in mind. In this work, we generalize to a more flexible framework where $\mathsf{F}$ and $\mathsf{G}$ can each be locally approximated at any point with a convex function.

A special case is the setting where $\mathsf{G}$ is convex while $\mathsf{F}$ is nonconvex, but $x \mapsto \mathsf{F}(Kx)$ is convex (i.e. if $\mathsf{F}$ is twice differentiable, then this is equivalent to assuming that $K^\top(\nabla^2\mathsf{F})K \succeq 0$ but $\nabla^2\mathsf{F} \not\succeq 0$), which may arise in many settings. In this case, the overall optimization problem, i.e. minimizing $\mathsf{F}(Kx) + \mathsf{G}(x)$, is a convex problem, that is, any local minimum is guaranteed to be a global minimum, and thus we might expect that this problem would be simple to optimize. Surprisingly, this may not be the case—if $\mathsf{F}$ is nondifferentiable, then we cannot use gradient descent on $\mathsf{F}$, while the nonconvexity of $\mathsf{F}$ means that existing primal/dual optimization techniques might not be applicable.

In the nonconvex, or more specifically, convex/concave setting, one of the most common techniques used in place of convex methods is the majorization/minimization approach [30, 16], where at each iteration, we work with a convex upper bound on the nonconvex objective function. Specifically, for the setting considered here, at each iteration we would choose some convex functions $\mathsf{F}^{(t)}, \mathsf{G}^{(t)}$ satisfying $\mathsf{F}^{(t)} \geq \mathsf{F}$ and $\mathsf{G}^{(t)} \geq \mathsf{G}$, then solve the modified optimization problem $\min_x\{\mathsf{F}^{(t)}(Kx) + \mathsf{G}^{(t)}(x)\}$. However, the modified optimization problem may itself be very challenging to solve in this setting, so we often cannot apply the majorization/minimization technique in a straightforward way. Our work combines the ideas of majorization/minimization with primal/dual techniques to handle this composite optimization problem.

## 1.1 The MOCCA algorithm

In this work, we propose the mirrored convex/concave (MOCCA) algorithm, which offers an approach that combines some of the techniques described above. We work with a primal/dual formulation of the problem, and incorporate a majorization/minimization-type step at each iteration. To motivate our method, we first present an existing method for the case where $\mathsf{F}$ and $\mathsf{G}$ are both convex, the Chambolle-Pock (CP) algorithm [7] (this method is closely related to other existing algorithms, which we discuss later on).

The CP algorithm is derived by considering the problem in a different form. From this point on, for clarity, we will use variables $x, y, z$ to denote points in $\mathbb{R}^d$ and $u, v, w$ to denote points in $\mathbb{R}^m$. Since we are considering the setting where $\mathsf{F}$ is convex, by duality we can write

$$\min_x\{\mathsf{F}(Kx) + \mathsf{G}(x)\} = \min_x \max_w\{\langle Kx, w\rangle - \mathsf{F}^*(w) + \mathsf{G}(x)\}\ .$$

where $\mathsf{F}^*$ is the conjugate to $\mathsf{F}$ [33], also known as the Legendre-Fenchel transform of $\mathsf{F}$, and is defined as

$$\mathsf{F}^*(w) = \max_v\{\langle w, v\rangle - \mathsf{F}(v)\}. \tag{2}$$

The primal variable $x$ and dual variable $w$ define a saddle-point problem. Given step sizes $\Sigma, \mathrm{T}$ which are positive diagonal matrices, the (preconditioned) CP algorithm [7, 32] iterates the following steps:

$$\begin{cases} x_{t+1} = \arg\min_x \left\{\langle Kx, w_t\rangle + \mathsf{G}(x) + \frac{1}{2}\left\|x - x_t\right\|^2_{\mathrm{T}^{-1}}\right\}, \\ w_{t+1} = \arg\min_w \left\{-\langle K\bar{x}_{t+1}, w\rangle + \mathsf{F}^*(w) + \frac{1}{2}\left\|w - w_t\right\|^2_{\Sigma^{-1}}\right\}, \end{cases} \tag{3}$$

where $\bar{x}_{t+1} = x_{t+1} + \theta(x_{t+1} - x_t)$ is an extrapolation term; here $\theta \geq 0$ (generally $\theta = 1$). Here the two norms are calculated via the definition $\|x\|_A := \sqrt{x^\top Ax}$ (for any positive semidefinite matrix $A \succeq 0$). When $\left\|\Sigma^{1/2}K\mathrm{T}^{1/2}\right\| < 1$, convergence properties for this algorithm have been proved, e.g. [7, 32, 14].[1] Setting $\theta = 0$

[1]The original form of the Chambolle-Pock algorithm [7], without preconditioning, can be obtained from (3) by replacing $\Sigma, \mathrm{T}$ with $\sigma\mathbf{I}, \tau\mathbf{I}$ for scalar step size parameters $\sigma, \tau > 0$, with convergence results proved if $\sigma\tau\left\|K\right\|^2 < 1$; however, in general the preconditioned form gives better performance and we only consider the preconditioned version here.

reduces to an earlier approach, the Primal-Dual Hybrid Gradient algorithm [11]; with $\theta = 1$, the CP algorithm is equivalent to a modification of ADMM (discussed later in Section 3.1).

We now ask how we could modify the discussed methods to handle nonconvexity of $\mathsf{F}$ and/or $\mathsf{G}$. One approach would be to approximate the problem with a convex optimization problem, that is, to take some approximation to $\mathsf{F}$ and $\mathsf{G}$ that are chosen to be a convex function. Of course, in general, there may not be a convex function that will provide a globally accurate approximation to a nonconvex function, but local convex approximations may be possible.

**Local convex approximations** Consider a family of approximations to the functions $\mathsf{F}$ and $\mathsf{G}$, indexed by $(z, v) \in \mathbb{R}^d \times \mathbb{R}^m$, the (primal) points at which the local approximations are taken. We write $\mathsf{F}_v : \mathbb{R}^m \to \mathbb{R}$ and $\mathsf{G}_z : \mathbb{R}^d \to \mathbb{R}$ for these approximations, and from this point on, we implicitly assume that the approximations satisfy:

$$\begin{cases} \mathsf{dom}(\mathsf{F}_v) = \mathsf{dom}(\mathsf{F}) \text{ and } \mathsf{dom}(\mathsf{G}_z) = \mathsf{dom}(\mathsf{G}) \text{ for any expansion points } (z,v) \in \mathsf{dom}(\mathsf{G}) \times \mathsf{dom}(\mathsf{F}); \\ \mathsf{F}_v \text{ and } \mathsf{G}_z \text{ are convex and continuous functions at any expansion points } (z,v) \in \mathsf{dom}(\mathsf{G}) \times \mathsf{dom}(\mathsf{F}); \\ \mathsf{F}_v \text{ and } \mathsf{G}_z \text{ are accurate up to first order, meaning that } \mathsf{F} - \mathsf{F}_v \text{ and } \mathsf{G} - \mathsf{G}_z \text{ are both differentiable,} \\ \qquad \text{and satisfy } \nabla(\mathsf{F} - \mathsf{F}_v)(v) = 0, \ \nabla(\mathsf{G} - \mathsf{G}_z)(z) = 0, \ \text{for any } (z,v) \in \mathsf{dom}(\mathsf{G}) \times \mathsf{dom}(\mathsf{F}). \end{cases} \tag{4}$$

(Here $\mathsf{dom}(\cdot)$ denotes the domain of a function, i.e. all points at which the function takes a finite value.) In particular this assumption implicitly requires that $\mathsf{F}$ and $\mathsf{G}$ both have convex domains.

As a special case, in some settings we can consider decomposing each function into the sum of a convex and a differentiable term, $\mathsf{F} = \mathsf{F}_{\mathsf{cvx}} + \mathsf{F}_{\mathsf{diff}}$ and $\mathsf{G} = \mathsf{G}_{\mathsf{cvx}} + \mathsf{G}_{\mathsf{diff}}$, as mentioned before; we can then take linear approximations to $\mathsf{F}_{\mathsf{diff}}$ and to $\mathsf{G}_{\mathsf{diff}}$ at the points $v$ and $z$, respectively, to obtain

$$\begin{cases} \mathsf{F}_v(w) := \mathsf{F}_{\mathsf{cvx}}(w) + \left[\mathsf{F}_{\mathsf{diff}}(v) + \langle w - v, \nabla \mathsf{F}_{\mathsf{diff}}(v)\rangle\right], \\ \mathsf{G}_z(x) := \mathsf{G}_{\mathsf{cvx}}(x) + \left[\mathsf{G}_{\mathsf{diff}}(z) + \langle x - z, \nabla \mathsf{G}_{\mathsf{diff}}(z)\rangle\right] . \end{cases} \tag{5}$$

In particular, if $\mathsf{F}_{\mathsf{diff}}$ and $\mathsf{G}_{\mathsf{diff}}$ are both concave (and thus $\mathsf{F}$ and $\mathsf{G}$ are each convex/concave), these two approximations are standard convex majorizations to $\mathsf{F}$ and $\mathsf{G}$ at the current estimates $w$ and $x$ (as might be used in a majorization/minimization algorithm in some settings).

**The algorithm** Since $\mathsf{F}_v, \mathsf{G}_z$ are convex, we can substitute them in place of $\mathsf{F}, \mathsf{G}$ in the iterations of the CP algorithm (3):

$$\begin{cases} x_{t+1} = \arg\min_x \left\{ \langle Kx, w_t\rangle + \mathsf{G}_z(x) + \frac{1}{2}\left\|x - x_t\right\|^2_{\mathrm{T}^{-1}} \right\}, \\ w_{t+1} = \arg\min_w \left\{ -\langle K\bar{x}_{t+1}, w\rangle + \mathsf{F}^*_v(w) + \frac{1}{2}\left\|w - w_t\right\|^2_{\Sigma^{-1}} \right\}. \end{cases} \tag{6}$$

We will see later on (in Section 3.3) that this formulation is closely related to the ADMM algorithm, and in fact in the special case given in (5), if $\mathsf{F}_{\mathsf{diff}}$ and $\mathsf{G}_{\mathsf{diff}}$ are concave, then MOCCA can be viewed as a special case of Bregman ADMM [41, 39].

Of course, a key question remains: which primal points $(z, v)$ should we use for constructing the convex approximations to $\mathsf{F}$ and to $\mathsf{G}$, at iteration $t$ of the algorithm? We find that, before solving for $(x_{t+1}, w_{t+1})$, we should use expansion points

$$(z, v) = \left(x_t, \Sigma^{-1}(w_{t-1} - w_t) + K\bar{x}_t\right) .$$

We will return shortly to the question of how these values were chosen; with this choice in place, the MOCCA algorithm is defined in Algorithm 1.

For reasons of stability, it may sometimes be desirable to update the expansion points $(z, v)$ only periodically, and we will incorporate this option into a more general version of MOCCA, given in Algorithm 2. Specifically, at the $t$th stage, we repeat the $(x, w)$ updates $L_t$ many times; we refer to these repeated updates as the "inner loop". In fact, the $t$th "inner loop" is simply running the CP algorithm for the convex problem $\min_x\{\mathsf{F}_{v_t}(Kx) + \mathsf{G}_{z_t}(x)\}$.

**Algorithm 1** MOCCA algorithm

**Input:** Functions $\mathsf{F}, \mathsf{G}$ with local convex approximations $\mathsf{F}_v, \mathsf{G}_z$, linear operator $K$, positive diagonal step size matrices $\Sigma, \mathrm{T}$, extrapolation parameter $\theta \in [0, 1]$.
**Initialize:** primal point $x_0 \in \mathbb{R}^d$, dual point $w_0 \in \mathbb{R}^m$, expansion points $(z_0, v_0) \in \mathbb{R}^d \times \mathbb{R}^m$.
**for** $t = 0, 1, 2, \ldots$ **do**

Update $x$ and $w$ variables: writing $\bar{x}_{t+1} = x_{t+1} + \theta(x_{t+1} - x_t)$, define

$$\begin{cases} x_{t+1} = \arg\min_x \left\{ \langle Kx, w_t \rangle + \mathsf{G}_{z_t}(x) + \frac{1}{2} \|x - x_t\|^2_{\mathrm{T}^{-1}} \right\}, \\ w_{t+1} = \arg\min_w \left\{ -\langle K\bar{x}_{t+1}, w \rangle + \mathsf{F}^*_{v_t}(w) + \frac{1}{2} \|w - w_t\|^2_{\Sigma^{-1}} \right\}. \end{cases}$$

Update expansion points: define

$$\begin{cases} z_{t+1} = x_{t+1}, \\ v_{t+1} = \Sigma^{-1}(w_t - w_{t+1}) + K\bar{x}_{t+1} \in \partial \mathsf{F}^*_{v_t}(w_{t+1}). \end{cases}$$

**until** some convergence criterion is reached.

Then, we average over the inner loop and calculate a single update of the expansion points $(z, v)$. Observe that, by setting $L_t = 1$, we do only a single $(x, w)$ update in each "inner loop" and thus have reduced to the basic form of MOCCA. In practice, the basic form (Algorithm 1) performs well, and the more stable version (Algorithm 2) is primarily proposed here for theoretical purposes, as some of our convergence guarantees do not hold for the basic version. However, in some settings the added stability does help empirically.

**Constraints and infinite function values** We remark that $\mathsf{F}, \mathsf{G}$ (and their approximations $\mathsf{F}_v, \mathsf{G}_z$) are allowed to take the value $+\infty$, for instance, to reflect a constraint. For example we might have $\mathsf{G}(x) = \delta\{\|x\|_2 \leq 1\}$, the convex indicator function taking the value $+\infty$ if the constraint $\|x\|_2 \leq 1$ is violated; this has the effect of imposing a constraint on the $x$ update step of our algorithm. Furthermore, in settings where $\mathsf{F}_v$ may not be strongly convex, its conjugate $\mathsf{F}^*_v$ may not be finitely valued; we would have $\mathsf{F}^*_v(w) = +\infty$ (for some, but not all, $w$). For instance if $\mathsf{F}_v(w) = \|w\|_1$ then $\mathsf{F}^*_v = \delta\{\|w\|_\infty \leq 1\}$, which has the effect of imposing a constraint on $w$ in the $w$ update step of our algorithm. Our theoretical results in this paper hold across all these settings, i.e. we do not assume that any of the functions $\mathsf{F}, \mathsf{G}, \mathsf{F}_v, \mathsf{G}_z, \mathsf{F}^*_v$ are everywhere finite, but instead work in the domains of these functions.

### 1.1.1 Simple special cases

Before discussing the implementation and behavior of MOCCA for general nonconvex/nonsmooth problems, we pause to illustrate that MOCCA can be viewed as a generalization of many existing techniques.

- If $\mathsf{F}, \mathsf{G}$ are both convex with easy proximal maps, then we can of course choose the trivial convex families $\mathsf{F}_v = \mathsf{F}$ and $\mathsf{G}_z = \mathsf{G}$; the MOCCA algorithm then reduces to the Chambolle-Pock [7] or Primal-Dual Hybrid Gradient [11] algorithm (depending on our choice of the extrapolation parameter $\theta$). These methods can handle composite objective functions with convex terms; MOCCA extends these methods to a setting with nonconvex terms.

- In the setting where we want to minimize a function $\mathsf{G}(x)$ which is a sum of a convex term and a differentiable term, $\mathsf{G}(x) = g(x) + h(x)$, we can show that MOCCA reduces to proximal gradient descent [27, 2] as a special case. Specifically, we define the approximations $\mathsf{G}_z(x) = g(x) + h(z) + \langle \nabla h(z), x - z \rangle$. In this setting there is no $\mathsf{F}$ function, and hence no $w$ variable; taking $\mathrm{T} = \tau \cdot \mathbf{I}_d$, the steps of Algorithm 1

**Algorithm 2** MOCCA algorithm (stable version with "inner loop")

**Input / Initialize:** Same as for Algorithm 1, along with inner loop lengths $L_1, L_2, \ldots$
**for** $t = 0, 1, 2, \ldots$ **do**
Define $(x_{t+1;0}, w_{t+1;0}) = (x_t, w_t)$.
**for** $\ell = 1, 2, \ldots, L_{t+1}$ **do**
Update $x$ and $w$ variables: writing $\bar{x}_{t+1;\ell} = x_{t+1;\ell} + \theta(x_{t+1;\ell} - x_{t+1;\ell-1})$, define

$$\begin{cases} x_{t+1;\ell} = \arg\min_x \left\{ \langle Kx, w_{t+1;\ell-1}\rangle + \mathsf{G}_{z_t}(x) + \frac{1}{2}\left\|x - x_{t+1;\ell-1}\right\|^2_{\mathrm{T}^{-1}} \right\}, \\ w_{t+1;\ell} = \arg\min_w \left\{ -\langle K\bar{x}_{t+1;\ell}, w\rangle + \mathsf{F}^*_{v_t}(w) + \frac{1}{2}\left\|w - w_{t+1;\ell-1}\right\|^2_{\Sigma^{-1}} \right\}. \end{cases}$$

**end for**
Define $(x_{t+1}, w_{t+1}) = \frac{1}{L_{t+1}} \sum_{\ell=1}^{L_{t+1}} (x_{t+1;\ell}, w_{t+1;\ell})$.
Update expansion points: define

$$\begin{cases} z_{t+1} = \frac{1}{L_{t+1}} \sum_{\ell=1}^{L_{t+1}} x_{t+1;\ell}, \\ v_{t+1} = \frac{1}{L_{t+1}} \sum_{\ell=1}^{L_{t+1}} \left(\Sigma^{-1}(w_{t+1;\ell-1} - w_{t+1;\ell}) + K\bar{x}_{t+1;\ell}\right). \end{cases}$$

**until** some convergence criterion is reached.

become

$$\begin{cases} x_{t+1} = \arg\min_x \left\{ \mathsf{G}_{z_t}(x) + \frac{1}{2\tau}\left\|x - x_t\right\|_2^2 \right\} = \arg\min_x \left\{ g(x) + \langle \nabla h(z_t), x\rangle + \frac{1}{2\tau}\left\|x - x_t\right\|_2^2 \right\}, \\ z_{t+1} = x_{t+1}, \end{cases}$$

which simplifies to the update scheme

$$x_{t+1} = \mathsf{Prox}_{\tau\cdot g}\left(x_t - \tau \cdot \nabla h(x_t)\right). \tag{7}$$

This is exactly the proximal gradient descent algorithm with step size $\tau$.

Proximal gradient descent can handle a function which combines a differentiable term with a convex term as long as the convex term has an easy proximal map; MOCCA extends this method to a setting where the convex terms lack an easy proximal map due to linear transformations, leading to composite optimization problems.

We will discuss other existing methods in more detail later on in Section 3.

#### 1.1.2 Step size parameters $\Sigma$ and $\mathrm{T}$

We now turn to the question of choosing the diagonal step size matrices $\Sigma$ and $\mathrm{T}$. As we will see later on, good convergence properties are attained when $\Sigma, \mathrm{T}$ are chosen sufficiently small, to satisfy $\left\|\Sigma^{1/2} K \mathrm{T}^{1/2}\right\| < 1$—this condition on $\Sigma$ and $\mathrm{T}$ is derived by [32] for the preconditioned CP algorithm, and appears in our theory as well. Here $\|\cdot\|$ is the matrix operator norm (i.e. the largest eigenvalue). To choose matrices that satisfy this requirement, [32] propose a parametrized family of choices: after fixing some parameter $\lambda > 0$, define[2]

$$\Sigma_{ii} = \frac{\lambda}{\sum_j |K_{ij}|} \quad \text{and} \quad \mathrm{T}_{jj} = \frac{\lambda^{-1}}{\sum_i |K_{ij}|}. \tag{8}$$

Empirically, we find that higher values of $\lambda$ are more stable but lead to slower convergence; it seems that the best choice is the smallest possible $\lambda$ such that the algorithm does not diverge. It may also be interesting to consider varying $\lambda$ adaptively over the iterations of the algorithm, but we do not study this extension here.

[2] In fact these choices for $\Sigma, \mathrm{T}$ satisfy the matrix norm constraint more weakly, with $\left\|\Sigma^{1/2} K \mathrm{T}^{1/2}\right\| \leq 1$ rather than a strict inequality, but this is sufficient in practice.

### 1.1.3 Understanding the choice of expansion points

We now return to our choice of the expansion points $(z, v) \in \mathbb{R}^d \times \mathbb{R}^m$. We will give an intuition for the choices of these points in the MOCCA algorithm. To examine this question, first consider the goal for optimization: we would like to find

$$x^\star \in \arg\min_x \{\mathsf{F}(Kx) + \mathsf{G}(x)\} \ ,$$

or if this problem is nonconvex then we may be satisfied to let $x^\star$ be a local minimizer or critical point of this objective function. We then need to find primal points $z \in \mathbb{R}^d$ and $v \in \mathbb{R}^m$, such that replacing $\mathsf{F}$ with $\mathsf{F}_v$ and $\mathsf{G}$ with $\mathsf{G}_z$ still yields the same solution, i.e. so that

$$x^\star \in \arg\min_x \{\mathsf{F}_v(Kx) + \mathsf{G}_z(x)\} \ . \tag{9}$$

Examining the first-order optimality conditions for each of these problems, it follows that we should set $(z, v) = (x^\star, Kx^\star)$ to ensure that (9) holds.

Of course, $x^\star$ is unknown and so we cannot set $(z, v) = (x^\star, Kx^\star)$. Instead, a logical approach would be to set $(z_t, v_t) = (x_t, Kx_t)$, before solving for $(x_{t+1}, w_{t+1})$. Then, hopefully, as $x_t$ converges to $x^\star$ we will also have $(z_t, v_t)$ converging to $(Kx^\star, x^\star)$. However, in practice, we find that this approach does not always perform as well as expected. Specifically, the problem lies with the choice $v_t = Kx_t$, relative to the primal/dual structure of the algorithm.

To understand why, imagine that $\mathsf{F}$ and $\mathsf{G}$ are actually convex, but we nonetheless are taking local approximations $\mathsf{F}_v$ and $\mathsf{G}_z$ (which are also convex), perhaps for computational reasons. Then we would like our $x$ and $w$ update steps (6) to coincide with the updates (3) of the original CP algorithm. Examining the optimality conditions, this will occur when

$$\partial\mathsf{G}(x_{t+1}) = \partial\mathsf{G}_{z_t}(x_{t+1}) \quad \text{and} \quad \partial\mathsf{F}^*(w_{t+1}) = \partial\mathsf{F}^*_{v_t}(w_{t+1}) \ . \tag{10}$$

(Here we are ignoring issues of multivalued subdifferentials since our aim is only to give intuition.) Using the definitions of $\mathsf{G}$ and $\mathsf{G}_{z_t}$, for the $x$ step our requirement in (10) is equivalent to

$$\nabla(\mathsf{G} - \mathsf{G}_{z_t})(x_{t+1}) = 0,$$

which will certainly hold if $z_t = x_{t+1}$ by our assumption (4) on the expansions $\mathsf{G}_{z_t}$. Since we have not yet solved for $x_{t+1}$, we instead choose the expansion point $z_t = x_t$ for the function $\mathsf{G}$, as previously proposed.

For the $w$ step, our outcome will be different. Subgradients satisfy a duality property, namely $w \in \partial\mathsf{F}(u)$ if and only if $u \in \partial\mathsf{F}^*(w)$ for any convex function $\mathsf{F}$ and its conjugate $\mathsf{F}^*$ (recall (2)). The requirement $\partial\mathsf{F}^*(w_{t+1}) = \partial\mathsf{F}^*_{v_t}(w_{t+1})$ in (10) therefore yields $w_{t+1} \in \partial\mathsf{F}(\partial\mathsf{F}^*_{v_t}(w_{t+1}))$ by this duality property, and so we have

$$\begin{aligned} w_{t+1} &\in \partial\mathsf{F}(\partial\mathsf{F}^*_{v_t}(w_{t+1})) \\ &= \partial\mathsf{F}_{v_t}(\partial\mathsf{F}^*_{v_t}(w_{t+1})) + \nabla(\mathsf{F} - \mathsf{F}_{v_t})(\partial\mathsf{F}^*_{v_t}(w_{t+1})) \\ &= w_{t+1} + \nabla(\mathsf{F} - \mathsf{F}_{v_t})(\partial\mathsf{F}^*_{v_t}(w_{t+1})) \end{aligned}$$

where the last step again holds from the duality property of subgradients. So, we see that we would like

$$\nabla(\mathsf{F} - \mathsf{F}_{v_t})(\partial\mathsf{F}^*_{v_t}(w_{t+1})) = 0$$

which, according to our assumption (4) on the expansions $\mathsf{F}_{v_t}$, will hold if

$$v_t \in \partial\mathsf{F}^*_{v_t}(w_{t+1}) \ .$$

In other words, we would like $v_t$ to be the primal point that corresponds to the dual point $w_{t+1}$—that is, the primal point that *mirrors* the dual point $w_{t+1}$. Of course, this is not possible since we have not yet computed $w_{t+1}$, and furthermore $v_t$ appears on both sides of this equation. Instead, we take $v_t \in \partial\mathsf{F}^*_{v_{t-1}}(w_t)$. Looking at the first-order optimality conditions for the update step for $w_t$, we see that we can satisfy this expression by choosing

$$v_t = \Sigma^{-1}(w_{t-1} - w_t) + K\bar{x}_t \in \partial\mathsf{F}^*_{v_{t-1}}(w_t) \ .$$

In fact, we will see in Section 3.3 that this choice for $v_t$ is very logical in light of the connection between the CP algorithm and the Alternating Direction Method of Multipliers (ADMM) [4].

For the stable form of MOCCA, Algorithm 2, our choice for expansion points $(z, v)$ takes an average over each inner loop, which we will see gives sufficient stability for our convergence results to hold.

#### 1.1.4 Checking convergence

Here we give a simple way to check whether the basic MOCCA algorithm, Algorithm 1, is near a critical point (e.g. a local minimum). (We treat this question more formally, for the more general Algorithm 2, in our theoretical results later on.) Due to the first-order accuracy of the convex approximations to $\mathsf{F}$ and $\mathsf{G}$ as specified in (4), a critical point $x \in \mathbb{R}^d$ for the objective function $\mathsf{F}(Kx) + \mathsf{G}(x)$ is characterized by the first-order condition

$$0 \in K^\top \partial \mathsf{F}_{Kx}(Kx) + \partial \mathsf{G}_x(x) \; . \tag{11}$$

Equivalently, we can search for a dual variable $w \in \mathbb{R}^m$ such that

$$w \in \partial \mathsf{F}_{Kx}(Kx) \quad \text{and} \quad -K^\top w \in \partial \mathsf{G}_x(x).$$

We can expand this condition to include additional variables $(z, v) \in \mathbb{R}^d \times \mathbb{R}^m$:

$$\begin{cases} w \in \partial \mathsf{F}_v(Kx) \;\Leftrightarrow\; Kx \in \partial \mathsf{F}_v^*(w), \\ -K^\top w \in \partial \mathsf{G}_z(x), \\ z = x, \\ v = Kx. \end{cases}$$

To check whether these conditions hold approximately, we can take the following "optimality gap":

$$\mathsf{OptimalityGap}(x, w, z, v) = \|Kx - \partial \mathsf{F}_v^*(w)\|_2^2 + \left\|-K^\top w - \partial \mathsf{G}_z(x)\right\|_2^2 + \|z - x\|_2^2 + \|v - Kx\|_2^2 \; .$$

Here, if any of the subdifferentials are multivalued, we can interpret these norms as choosing some element of the corresponding subdifferentials. Now we consider the value of this gap at an iteration of the MOCCA algorithm (in its original form, Algorithm 1). By the definitions of $x_{t+1}, w_{t+1}, z_t, v_t$, we can show that

$$\begin{cases} 0 \in K^\top w_t + \partial \mathsf{G}_{z_t}(x_{t+1}) + \mathrm{T}^{-1}(x_{t+1} - x_t), \\ 0 \in -K\bar{x}_{t+1} + \partial \mathsf{F}_{v_t}^*(w_{t+1}) + \Sigma^{-1}(w_{t+1} - w_t), \\ z_t = x_t, \\ v_t = \Sigma^{-1}(w_t - w_{t-1}) + K\bar{x}_t. \end{cases}$$

Therefore, plugging these calculations in to the definition of the optimality gap, we see that

$$\begin{aligned} &\mathsf{OptimalityGap}(x_{t+1}, w_{t+1}, z_t, v_t) \\ &\quad = \left\|K(x_t - x_{t+1}) + \Sigma^{-1}(w_t - w_{t+1})\right\|_2^2 + \left\|-K(w_t - w_{t+1}) + \mathrm{T}^{-1}(x_t - x_{t+1})\right\|_2^2 \\ &\qquad + \|x_t - x_{t+1}\|_2^2 + \left\|K(x_{t-1} - 2x_t + x_{t+1}) + \Sigma^{-1}(w_{t-1} - w_t)\right\|_2^2 \\ &\qquad = \mathcal{O}\left( \left\| \begin{pmatrix} x_{t-1} - x_t \\ w_{t-1} - w_t \end{pmatrix} \right\|_2^2 + \left\| \begin{pmatrix} x_t - x_{t+1} \\ w_t - w_{t+1} \end{pmatrix} \right\|_2^2 \right) \; . \end{aligned}$$

In other words, if the change in $(x_t, w_t)$ converges to zero as $t \to \infty$, then the optimality gap is also converging to zero.

### 1.1.5 Preview of theoretical results

We present two theoretical results in this work. The first is fairly standard in the related literature: in Theorem 1 we show that if the algorithm does converge to a point, then we have reached a critical point of the original optimization problem. (Since the simple form, Algorithm 1, is a special case of the stable form, Algorithm 2, we prove this result for the stable algorithm only.)

The novelty of our theory lies in our convergence guarantee, given in Theorem 2, where we prove that the stable form of MOCCA, given in Algorithm 2, is guaranteed to converge to a nearly-globally-optimal solution, under some assumptions on convexity and curvature. Specifically, we consider a scenario where convexity and nonconvexity in $\mathsf{F}$ and $\mathsf{G}$ counterbalance each other, so that the overall function

$$x \mapsto \mathsf{F}(Kx) + \mathsf{G}(x) \tag{12}$$

is itself either strongly convex or satisfies restricted strong convexity assumptions (which we will discuss in detail in Section 4.2.1). It is important to note that even the globally convex setting is by no means trivial—even if (12) is strongly convex, if $\mathsf{F}$ itself is nonconvex it may be the case that ADMM and other primal/dual or alternating minimization algorithms diverge or converge to the wrong solution, as we discuss later in Section 3.2. Crucially, our results allow $\mathsf{F}$ and $\mathsf{G}$ to be nondifferentiable as well as nonconvex, a setting that is necessary in practice but is not covered by existing theory.

### 1.1.6 Outline of paper

The remainder of the paper is organized as follows. In Section 2, we present several important applications where the minimization problem considered here, with a nonconvex composite objective function as in (1), arises naturally: regression problems with errors in covariates, isotropic total variation penalties, nonconvex total variation penalties, and image reconstruction problems in computed tomography (CT) imaging. In Section 3 we give background on several types of existing algorithms for convex and nonconvex composite objective functions, and compare a range of existing results to the work presented here. Theoretical results on the convergence properties of our algorithm are given in Section 4. We study the empirical performance of MOCCA in Section 5. Proofs are given in Section 6, with technical details deferred to the Appendix. In Section 7 we discuss our findings and outline directions for future research.

# 2 Applications

We now highlight several applications of the MOCCA algorithm, in high-dimensional statistics and in imaging.

## 2.1 Regression with errors in variables

Recent work by Loh and Wainwright [22] considers optimization for nonconvex statistical problems, proving that under some special conditions, nonconvexity may not pose a challenge to recovering optimal parameters. In particular, they consider the following example [21, 22]: suppose that we observe a response $y \in \mathbb{R}^n$ which is generated with a Gaussian linear model,

$$b = Ax_{\mathsf{true}} + \epsilon \text{ with } \epsilon \sim N(0, \sigma^2 \mathbf{I}_t) \ ,$$

where $A \in \mathbb{R}^{n\times d}$ is a design matrix and $x_{\mathsf{true}} \in \mathbb{R}^d$ is the unknown vector of coefficients. In this case, we might seek to recover $x_{\mathsf{true}}$ with the least squares estimator, perhaps with some penalty added to promote some desired structure in $x$,

$$\widehat{x} = \arg\min_x \left\{ \frac{1}{2n} \|b - Ax\|_2^2 + \mathsf{Penalty}(x) \right\} = \arg\min_x \left\{ \frac{1}{2} x^\top \left( \frac{A^\top A}{n} \right) x - x^\top \left( \frac{A^\top b}{n} \right) + \mathsf{Penalty}(x) \right\} \ . \tag{13}$$

In some settings, however the design matrix $A$ itself may not be available with perfect accuracy. Instead, suppose we observe

$$Z = A + W$$

where $W \perp\!\!\!\perp A$ has independent mean-zero entries, with $\mathbb{E}\left[W_{ij}^2\right] = \sigma_A^2$ for all $i, j$. In this case, a naive approach might be to substitute $Z$ for $A$ in (13), before finding the minimizer. However, unless $\sigma_A^2$ is negligible, this may not produce a good approximation to $\widehat{x}$ since we have

$$\mathbb{E}\left[Z^\top Z \mid A\right] = \mathbb{E}\left[(A+W)^\top (A+W) \mid A\right] = A^\top A + \mathbb{E}\left[W^\top W\right] = A^\top A + n\sigma_A^2 \mathbf{I}_d\ .$$

In contrast, for the linear term, we have $\mathbb{E}\left[Z^\top b \mid A\right] = A^\top b$, as desired. To correct for the bias in $Z^\top Z$, we should take

$$\widehat{x}_{\mathsf{noisy}} = \arg\min_x \left\{\mathcal{L}(x) + \mathsf{Penalty}(x)\right\}\ ,$$

where

$$\mathcal{L}(x) := \frac{1}{2} x^\top \left(\frac{Z^\top Z}{n} - \sigma_A^2 \mathbf{I}_d\right) x - x^\top \left(\frac{Z^\top b}{n}\right)\ .$$

Of course, this optimization problem is no longer convex due to the negative quadratic term, and in particular, for a Lipschitz penalty and a high-dimensional setting ($n < d$), the value tends to $-\infty$ as $x$ grows large in any direction in the null space of $Z$. Remarkably, Loh and Wainwright [22] show that, if $x_{\mathsf{true}}$ is sparse and $\mathsf{Penalty}(x)$ is similar to the $\ell_1$ norm, then as long as $(Z^\top Z)$ satisfies a restricted strong convexity assumption (as is standard in the sparse regression literature), then $x_{\mathsf{true}}$ can be accurately recovered from *any* local minimum or critical point of the constrainted optimization problem

$$\widehat{x}_{\mathsf{noisy}} = \arg\min_{\|x\|_1 \le R} \left\{\mathcal{L}(x) + \mathsf{Penalty}(x)\right\}\ . \tag{14}$$

The approach taken by [22] is to perform proximal gradient descent, with steps taking the form

$$x_{t+1} = \arg\min_x \left\{ \frac{1}{2} \left\| x - \left(x_t - \frac{1}{\eta} \nabla \mathcal{L}(x_t)\right)\right\|_2^2 + \frac{1}{\eta} \mathsf{Penalty}(x)\right\}$$

where $\frac{1}{\eta}$ is a step size parameter. When $\mathsf{Penalty}(x)$ is (some multiple of) the $\ell_1$ norm, or some other function with a simple proximal operator (that is, it is simple to compute $\arg\min_x\{\frac{1}{2}\|x - z\|_2^2 + \mathsf{Penalty}(x)\}$ for any $z$), this algorithm is very efficient.

**Total variation and generalized convex $\ell_1$ penalties** Consider a setting where the penalty is a generalized $\ell_1$ norm,

$$\mathsf{Penalty}(x) = \nu \cdot \|Kx\|_1$$

for some matrix $K \in \mathbb{R}^{m \times d}$. In particular, if $K$ is the one-dimensional differences matrix $\nabla_{\mathrm{1d}} \in \mathbb{R}^{(p-1)\times p}$, which has entries $(\nabla_{\mathrm{1d}})_{ii} = 1$ and $(\nabla_{\mathrm{1d}})_{i,i+1} = -1$ for each $i$, then this defines the (one-dimensional) total variation norm on $x$, $\|x\|_{\mathsf{TV}} = \|\nabla_{\mathrm{1d}} x\|_1$; this method is also known as the fused Lasso [37]. We can also consider a two- or three-dimensional total variation norm, $K = \nabla_{\mathrm{2d}}$ or $K = \nabla_{\mathrm{3d}}$, defined analogously as the differences matrix for a two- or three-dimensional grid. Total variation type penalties are commonly used in imaging applications and many other fields to obtain solutions that are locally constant or locally smooth.

In this setting, proximal gradient descent is not practical except in some special cases, such as when $K$ is diagonal, because the proximal operator $\arg\min_x\{\frac{1}{2}\|z - x\|_2^2 + \nu \cdot \|Kx\|_1\}$ does not have a closed form solution for general $K$ and would itself require an iterative algorithm to be run to convergence. For a total variation penalty on a one-dimensional grid of points, e.g. $K = \nabla_{\mathrm{1d}}$, some fast algorithms do exist for the proximal map [17]. Additional methods for convex problems with two-dimensional total variation and related penalties such as total variation over a graph can be found in [6, 42, 43]. We are not aware of a non-iterative algorithm for general $K$. Here we apply MOCCA to allow for arbitrary $K$ and for a nonconvex loss term $\mathcal{L}(x)$.

**Applying the MOCCA algorithm** We consider applying the MOCCA algorithm to this nonconvex optimization problem with $\mathsf{Penalty}(x) = \nu \left\|Kx\right\|_1$. We define the convex function

$$\mathsf{F}(w) = \nu \left\|w\right\|_1$$

with the trivial convex approximations $\mathsf{F}_v(w) = \mathsf{F}(w)$ at any expansion point $v \in \mathbb{R}^m$. We also let

$$\mathsf{G}(x) = \begin{cases} \mathcal{L}(x), & \text{if } \left\|x\right\|_1 \leq R, \\ +\infty, & \text{if } \left\|x\right\|_1 > R, \end{cases}$$

with convex approximations given by taking the linear approximation to the loss,

$$\mathsf{G}_z(x) = \begin{cases} \mathcal{L}(z) + \langle x - z, \nabla \mathcal{L}(z) \rangle, & \text{if } \left\|x\right\|_1 \leq R, \\ +\infty, & \text{if } \left\|x\right\|_1 > R, \end{cases}$$

for any expansion point $z \in \mathbb{R}^d$. Then the optimization problem (14) can be expressed as

$$\widehat{x}_{\mathsf{noisy}} = \arg\min_x \left\{ \mathsf{F}(Kx) + \mathsf{G}(x) \right\} .$$

Applying Algorithm 1, the update steps take the form

$$\begin{cases} x_{t+1} = \arg\min_{\left\|x\right\|_1 \leq R} \left\{ \left\| x - \left[ x_t - \mathrm{T} \left( \nabla \mathcal{L}(x_t) + K^\top w_t \right) \right] \right\|^2_{\mathrm{T}^{-1}} \right\}, \\ w_{t+1} = \mathsf{Truncate}_\nu \left( w_t + \Sigma K \bar{x}_{t+1} \right), \\ z_{t+1} = x_{t+1}, \end{cases} \tag{15}$$

where $\mathsf{Truncate}_\nu(w)$ truncates the entries of the vector $w$ to the range $[-\nu, \nu]$. Note that the $x$ update step is a simple shrinkage step and therefore easy to solve (and $\nabla \mathcal{L}(x_t)$ is simple to compute), while the $w$ and $z$ updates are computationally trivial.

As a second option, we can incorporate more convexity into our approximations $\mathsf{G}_z$ by taking

$$\mathsf{G}_z(x) = \begin{cases} \mathcal{L}(x) + \frac{\sigma_A^2}{2} \left\|x - z\right\|_2^2, & \text{if } \left\|x\right\|_1 \leq R, \\ +\infty, & \text{if } \left\|x\right\|_1 > R, \end{cases}$$

which is convex since $\mathcal{L}(x)$ has negative curvature bounded by $\frac{\sigma_A^2}{2}$. In this case, after simplifying, our update steps become

$$\begin{cases} x_{t+1} = \arg\min_{\left\|x\right\|_1 \leq R} \left\{ \left\| x - \left[ x_t - \left( \mathrm{T}^{-1} + \frac{Z^\top Z}{n} \right)^{-1} \left( \nabla \mathcal{L}(x_t) + K^\top w_t \right) \right] \right\|^2_{\left( \mathrm{T}^{-1} + \frac{Z^\top Z}{n} \right)} \right\}, \\ w_{t+1} = \mathsf{Truncate}_\nu \left( w_t + \Sigma K \bar{x}_{t+1} \right), \\ z_{t+1} = x_{t+1}. \end{cases} \tag{16}$$

While the $x$ update step may appear difficult due to the combination of the non-diagonal matrix $\left( \mathrm{T}^{-1} + \frac{Z^\top Z}{n} \right)$ which scales the norm, combined with the $\ell_1$ constraint, in practice the constraint $R$ is chosen to be large so that it is inactive in all or most steps; the $x$ update step is then solved by $x_{t+1} = x_t - \left( \mathrm{T}^{-1} + \frac{Z^\top Z}{n} \right)^{-1} \left( \nabla \mathcal{L}(x_t) + K^\top w_t \right)$.

An important point is that MOCCA can be applied to this problem for arbitrary $K$, including a difference operator such as $\nabla_{\mathsf{2d}}$ or $\nabla_{\mathsf{3d}}$; in contrast, proximal gradient descent can only be performed approximately except for certain special cases, as mentioned above. We explore this setting's theoretical properties in Section 4.2.3, and give empirical results for this problem in Section 5.

## 2.2 Isotropic total variation and generalized $\ell_1/\ell_2$ penalties

For locally constant images or signals in two dimensions, the form of the total variation penalty given above is known as "anisotropic", meaning that it imposes a sparsity pattern which is specific to the alignment of the image onto a horizontal and vertical axis. In contrast, the isotropic total variation penalty [34], on an image $x$ parametrized with values $x_{i,j}$ at grid location $(i,j)$, is given by

$$\|x\|_{\mathsf{isoTV}} = \sum_{(i,j)} \sqrt{(x_{i,j} - x_{i,j+1})^2 + (x_{i,j} - x_{i+1,j})^2}.$$

Optimization methods for the denoising problem with an isotropic total variation penalty, i.e. problems of the form $\min_x \left\{\frac{1}{2}\|b - x\|_2^2 + \nu \cdot \|x\|_{\mathsf{isoTV}}\right\}$, were studied in [8]. In practice the isotropic penalty is often preferred as it leads to smoother contours, avoiding the artificial horizontal or vertical edges that may result from anistropic total variation regularization.

The isotropic total variation penalty can be generalized to penalties of the form

$$\mathsf{Penalty}(x) = \nu \cdot \sum_{\ell=1}^{L} \|K_\ell x\|_2 ,$$

where each $K_\ell \in \mathbb{R}^{m_\ell \times d}$ is some fixed matrix. To see how this determines an isotropic total variation penalty in two dimensions, we let $\ell$ index all locations $(i,j)$; then the corresponding matrix $K_\ell$ has two rows, which when multiplied by $x$, extract the differences $x_{i,j} - x_{i,j+1}$ and $x_{i,j} - x_{i+1,j}$. To see how this specializes to the usual generalized $\ell_1$ penalty, $\|Kx\|_1$ for some fixed matrix $K \in \mathbb{R}^{m\times d}$, simply take $K_\ell$ to be the $\ell$th row of the matrix $K$ for each $\ell = 1, \dots, L = m$.

**Applying the MOCCA algorithm** We now show the steps of the MOCCA algorithm to the problem of minimizing

$$\mathcal{L}(x) + \nu \cdot \sum_{\ell=1}^{L} \|K_\ell x\|_2 ,$$

where $\mathcal{L}(x)$ is a differentiable likelihood term (such the nonconvex likelihood for regression with errors in variables as above). We define the convex function

$$\mathsf{F}(w) = \nu \cdot \sum_{\ell}^{L} \|w_{B_\ell}\|_2 ,$$

where $w_{B_\ell}$ is understood to be the $\ell$th block of the vector $w$, i.e. $w_{B_\ell} = (w_{m_1+\dots+m_{\ell-1}+1}, \dots, w_{m_1+\dots+m_\ell})$. We take the trivial convex approximations $\mathsf{F}_v(w) = \mathsf{F}(w)$ at any expansion point $v \in \mathbb{R}^m$, and define $\mathsf{G}(x)$ and $\mathsf{G}_z(x)$ exactly as before (as for the previous application, we allow the option of restricting to $\|x\|_1 \leq R$ if desired). Applying Algorithm 1, the update steps take the form

$$\begin{cases} x_{t+1} = \arg\min_{\|x\|_1 \leq R} \left\{ \left\| x - \left[x_t - \mathrm{T}\left(\nabla\mathcal{L}(x_t) + K^\top w_t\right)\right]\right\|^2_{\mathrm{T}^{-1}} \right\}, \\ w_{t+1} = \mathsf{Truncate}_{\nu\cdot(\mathbb{B}_{m_1}\times\dots\times\mathbb{B}_{m_L})}\left(w_t + \Sigma K \bar{x}_{t+1}\right), \\ z_{t+1} = x_{t+1}, \end{cases}$$

where $\mathsf{Truncate}_{\nu\cdot(\mathbb{B}_{m_1}\times\dots\times\mathbb{B}_{m_L})}(w)$ projects each block $w_{B_\ell} \in \mathbb{R}^{m_\ell}$ of the vector $w$ to the ball of radius $\nu$, $\nu \cdot \mathbb{B}_{m_\ell} \subseteq \mathbb{R}^{m_\ell}$ (here $\mathbb{B}_{m_\ell}$ is the unit ball of dimension $m_\ell$ in the $\ell_2$ norm). Specifically, the $w$ update step can be computed for each block as

$$(w_{t+1})_{B_\ell} = (w_t + \Sigma K \bar{x}_{t+1})_{B_\ell} \cdot \min\left\{1, \frac{\nu}{\left\|(w_t + \Sigma K\bar{x}_{t+1})_{B_\ell}\right\|_2}\right\}$$

for $\ell = 1, \dots, L$, and is therefore trivial to compute.

## 2.3 Nonconvex total variation penalties

As discussed in Section 2.1, total variation penalties are common in many applications where the underlying signal exhibits smooth or locally constant spatial structure (in one, two, or more dimensions). In convex optimization, we are faced with a well-known tradeoff between sparsity and bias—using $\nu \cdot \|x\|_{\mathsf{TV}}$ as our penalty function for some parameter $\nu > 0$, we want to be sure to choose $\nu$ large enough that the resulting solution is total-variation-sparse, to avoid overfitting when the sample size is small relative to the dimension of the image; however, larger $\nu$ leads to increased shrinkage of the signal, leading to an estimate that is biased towards zero. One way to avoid this tradeoff is to use a nonconvex penalty, which should behave like the total variation norm in terms of promoting sparsity, but reduce the amount of shrinkage for larger signal strength. In this section, we will use $\nabla_{\text{TV}}$ to denote the differences matrix in the appropriate space (e.g. $\nabla_{\text{TV}} = \nabla_{\text{2d}}$ in two dimensions), so that $\|x\|_{\mathsf{TV}} = \|\nabla_{\text{TV}} x\|_1$.

For sparse regression problems (i.e. where the signal $x$ is itself sparse, rather than sparsity in $\nabla_{\text{TV}} x$), many nonconvex alternatives to the $\ell_1$ norm penalty $\|x\|_1$ have been studied, demonstrating more accurate signal recovery empirically as well as theoretical properties of reduced bias, such as the SCAD penalty [12], the $\ell_q$ penalty (penalizing $\sum_i |x_i|^q$ for some $q \in (0,1)$) [18, 9], and the MCP penalty which seeks to minimize concavity while avoiding bias [44]. Another option is to use a reweighted $\ell_1$ norm [5] where signals estimated to be large at the first pass, are penalized less in the next pass to reduce bias in their estimates; in fact, Candès et al. [5] shows that this procedure is related to a nonconvex log-sum penalty, given by penalizing each component of $x_i$ as $\log(|x_i| + \epsilon)$ for some fixed $\epsilon > 0$. For the problem of total variation sparsity, a variety of nonconvex approaches have also been studied, including applying SCAD [10], an $\ell_q$ norm penalty for $0 < q < 1$ [36, 23], or a log-sum total variation penalty [35, 31] to the total variation sparsity setting.

We now consider the problem applying a log-sum penalty to the problem of total variation sparsity. Here we consider the form of this penalty is given by

$$\mathsf{logTV}_\beta(x) = \mathsf{logL1}_\beta(\nabla_{\text{TV}} x) \text{ where } \mathsf{logL1}_\beta(w) = \sum_i \beta \log\left(1 + |w_i|/\beta\right)$$

where $\beta > 0$ is a nonconvexity parameter. (We can also consider applying this nonconvex penalty to the isotropic version of total variation, as discussed in Section 2.2, but for simplicity we do not give that version explicitly here.)

To understand this function, observe that for any $t$, the function

$$t \mapsto \beta \log(1 + |t|/\beta)$$

is approximately equal to $|t|$ when $t \approx 0$ (that is, near zero it behaves like the $\ell_1$ norm), but is nonconvex and penalizes large values of $t$ much less heavily than an absolute value penalty of $|t|$. The parameter $\beta$ controls the amount of nonconvexity; small $\beta$ gives a highly nonconvex penalty, while for large $\beta$ the penalty is nearly convex, with $\mathsf{logTV}_\beta(x) \approx \|x\|_{\mathsf{TV}}$; see Figure 1 for an illustration.

Consider the problem of minimizing an objective function

$$\mathcal{L}(x) + \nu \cdot \mathsf{logTV}_\beta(x) \ , \tag{17}$$

where $\mathcal{L}(x)$ is some likelihood or loss term. In the image denoising setting, where $\mathcal{L}(x) = \frac{1}{2}\|y - x\|_2^2$ (i.e. when $y$ is a noisy observation of the signal $x$), Parekh and Selesnick [31] approach this optimization problem with a majorization/minimization algorithm, iterating the steps: (1) find a majorization of $\mathsf{logTV}_\beta(x)$ at the current estimate $x_t$, which takes the form of a reweighted TV norm, (2) compute $x_{t+1}$ as the minimizer of $\mathcal{L}(x) + \nu \cdot$ (the majorized penalty). In other settings, however, for a general loss $\mathcal{L}(x)$, step (2) may not be possible.

We now show how the MOCCA algorithm can be used to optimize objective functions of the form (17). For simplicity we show the steps for the case that $\mathcal{L}(x)$ is convex and has a simple proximal operator, but this can be generalized as needed.

First, we define a new function

$$h_\beta(w) = \mathsf{logL1}_\beta(w) - \|w\|_1 = \sum_i \beta \log\left(1 + |w_i|/\beta\right) - \|w\|_1 \ .$$

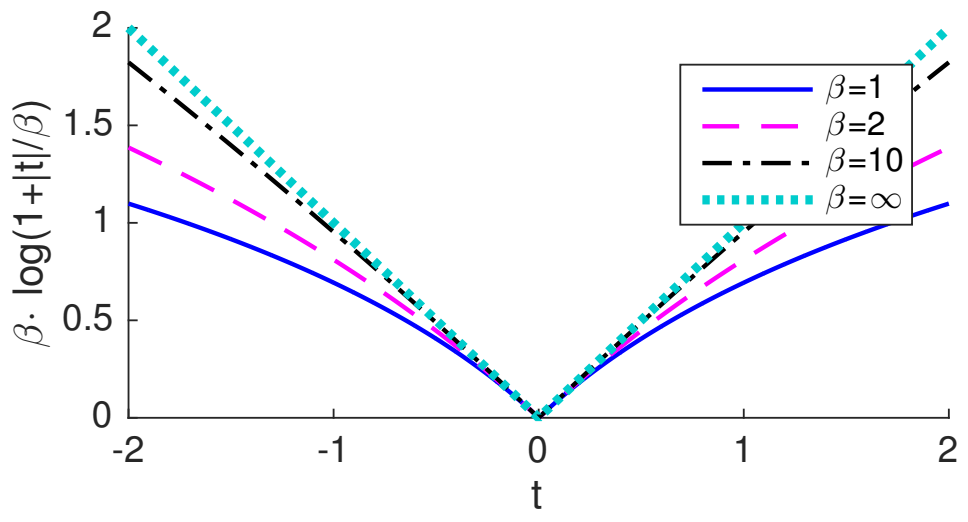

Figure 1: Illustration of the nonconvex sparsity-promoting penalty discussed in Section 2.3. The figure plots the function $t \mapsto \beta \cdot \log(1 + |t|/\beta)$ across a range of values of $t$, for $\beta \in \{1, 2, 10, \infty\}$; for $\beta = \infty$, we should interpret this as the absolute value function, $t \mapsto |t|$. We see that all the functions appear similar for $t \approx 0$, with a nondifferentiable point at $t = 0$ which ensures sparsity when this function is used as a penalty. For larger values of $t$, smaller values of $\beta$ correspond to greater nonconvexity.

We note an important property of this function: $h_\beta(w)$ is differentiable with $(\nabla h_\beta(w))_i = -\frac{w_i}{\beta+|w_i|}$. For a first approach, we will take

$$\mathsf{F}(w) = \nu \cdot \mathsf{logL1}_\beta(w) \text{ and } \mathsf{G}(x) = \mathcal{L}(x).$$

Now we define a local convex approximation to $\mathsf{F}$ by writing $\mathsf{F} = \nu\left(\|\cdot\|_1 + h_\beta\right)$ and taking the linear approximation to $h_\beta$, namely,

$$\mathsf{F}_v(w) = \nu \cdot \|w\|_1 + \nu\left[h_\beta(v) + \langle w - v, \nabla h_\beta(v)\rangle\right].$$

And, since by assumption $\mathcal{L}(x)$ is convex and has a simple proximal operator, we can simply define

$$\mathsf{G}_z(x) = \mathsf{G}(x) = \mathcal{L}(x).$$

Then the objective function (17) is equal to $\mathsf{F}(\nabla_{\mathrm{TV}}x) + \mathsf{G}(x)$, and can be minimized with MOCCA. Applying Algorithm 1, the update steps take the form

$$\begin{cases} x_{t+1} = \arg\min_x \left\{\mathcal{L}(x) + \frac{1}{2}\left\|x - \left[x_t - \mathrm{T}\nabla_{\mathrm{TV}}^\top w_t\right]\right\|^2_{\mathrm{T}^{-1}}\right\}, \\ w_{t+1} = \mathsf{Truncate}_\nu\left(w_t + \Sigma\nabla_{\mathrm{TV}}\bar{x}_{t+1} - \nu\nabla h_\beta(v_t)\right) + \nu\nabla h_\beta(v_t), \\ v_{t+1} = \Sigma^{-1}(w_t - w_{t+1}) + K\bar{x}_{t+1}. \end{cases}$$

Of course, we have the flexibility to arrange the functions differently if we wish—for example, we could instead define

$$\mathsf{F}(w) = \nu \cdot \|w\|_1 \text{ and } \mathsf{G}(x) = \mathcal{L}(x) + \nu \cdot h_\beta(\nabla_{\mathrm{TV}}x).$$

In this case, we will define the local approximations as

$$\mathsf{F}_v(w) = \mathsf{F}(w) = \nu \cdot \|w\|_1$$

and

$$\mathsf{G}_z(w) = \mathcal{L}(x) + \nu\left[h_\beta(\nabla_{\mathrm{TV}}z) + \langle\nabla_{\mathrm{TV}}(x - z), \nabla h_\beta(\nabla_{\mathrm{TV}}z)\rangle\right].$$

In this case the objective function (17) is equal to $\mathsf{F}(\nabla_{\mathrm{TV}}x) + \mathsf{G}(x)$ as before. In this case, the update steps are

$$\begin{cases} x_{t+1} = \arg\min_x \left\{\mathcal{L}(x) + \frac{1}{2}\left\|x - \left[x_t - \mathrm{T}(\nu\nabla h_\beta(\nabla_{\mathrm{TV}}z_t) + \nabla_{\mathrm{TV}}^\top w_t)\right]\right\|^2_{\mathrm{T}^{-1}}\right\}, \\ w_{t+1} = \mathsf{Truncate}_\nu\left(w_t + \Sigma\nabla_{\mathrm{TV}}\bar{x}_{t+1}\right), \\ z_{t+1} = x_{t+1}. \end{cases}$$

However, in a sense this decomposition is less natural as it splits the penalty $\mathsf{logTV}_\beta(x)$ across $\mathsf{F}$ and $\mathsf{G}$, and in fact in our experiments in Section 5, we will see that this second formulation gives poorer convergence results when the MOCCA algorithm is applied.

## 2.4 Application to CT image reconstruction

The initial motivation for developing the algorithm presented here, is a problem arising in computed tomography (CT) imaging. Here we briefly summarize the problem and our approach via the MOCCA algorithm; this setting, and detailed results, are described more fully in our forthcoming work [1].

In CT imaging, an X-ray beam is sent along many rays through the object of interest. Typically, the measurement is the total energy that has passed through the object, along each ray; comparing the energy retrieved against the energy of the entering X-ray beam, gives information about the materials inside the object, since different materials have different beam attenuation properties. A recent technological development is the photon counting detector, which measures the raw number of photons that successfully pass through the object rather than the total integrated energy. As a first pass, this transmission model can be written as follows, where the location $\vec{r}_\ell(t)$ inside the object parametrizes the ray $\ell$:

$$\mathsf{Count}_\ell \sim \mathsf{Poisson}\left(\int_{\text{energy } E} \begin{pmatrix}\text{beam}\\ \text{intensity at}\\ \text{energy } E\end{pmatrix} \cdot \begin{pmatrix}\text{detector}\\ \text{sensitivity}\\ \text{to energy } E\end{pmatrix} \cdot \exp\left\{-\int_t \mu(E, \vec{r}_\ell(t))\ \mathsf{d}t\right\} \mathsf{d}E\right), \quad (18)$$

where the only unknowns are the coefficients $\mu = (\mu(E, \vec{r}))$, indexed over energy level $E$ and location $\vec{r}$ inside the object, where $\mu(E, \vec{r})$ is the attenuation for photons at energy $E$ at the location $\vec{r}$. Higher values of $\mu(E, \vec{r})$ indicate that a photon is more likely to be absorbed; the expression $\exp\{-\int_t \mu(E, \vec{r}_\ell(t))\ \mathsf{d}t\}$ determines, for a photon at energy $E$ that enters the object along the trajectory defined by ray $\ell$, the probability that the photon will not be absorbed by the object (i.e. will pass through the object). These coefficients $\mu(E, \vec{r})$ can be further decomposed as

$$\mu(E, \vec{r}) = \sum_{\text{materials } m} \mu_m(E) \cdot x_m(\vec{r}) \quad (19)$$

where $\mu_m(E)$ is a known quantity determining the absorption properties of material $m$ at energy $E$, while $x_m(\vec{r})$ is unknown, representing the amount of material $m$ that is present at location $\vec{r}$. In practice, the object space is discretized into pixels, so that $x$ is finite-dimensional.

In this application, the optimization problem is then given by

$$\widehat{x} = \arg\min_x \left\{\mathcal{L}(x) + \sum_{\text{materials } m} (\text{total variation constraint on material map } x_m)\right\},$$

perhaps with other constraints added as well (e.g. nonnegativity of the material maps). Here $\mathcal{L}(x)$ is the negative log-likelihood of $x$ given the Poisson model for the observed photon counts as a function of $x$, given by (18) and (19). $\mathcal{L}(x)$ is a nonconvex function due to the integration across energy levels. Both $\mathcal{L}(x)$ and the total variation constraints are better represented as functions of linear transformations of $x$; the presence of these multiple terms, including nonconvexity (from $\mathcal{L}$) and nondifferentiability (from total variation), mean that existing methods cannot be applied to solve this optimization problems.

The MOCCA algorithm, applied to this problem, gives strong performance in terms of fast convergence and accurate image reconstruction on simulated and real imaging data. We do not give the details of the algorithm implementation or any empirical results in this paper, but instead refer the reader to our work in [1] for more details on the method and for empirical results on simulated CT image data (results on real data are forthcoming).

# 3 Background: optimizing convex and nonconvex composite functions

In this sections, we give background on several related algorithms for solving the minimization problem in the simpler setting where $\mathsf{F}$ and $\mathsf{G}$ are convex, and in the more challenging nonconvex setting.

### 3.1 Optimization when $\mathsf{F}$ is convex

When $\mathsf{F}$ is convex, a variety of existing methods are available to recover the (possibly non-unique) minimizer

$$x^\star = \arg\min\{\mathsf{F}(Kx) + \mathsf{G}(x)\} \ .$$

In many settings, the functions $\mathsf{F}$ and $\mathsf{G}$ may each have easily computable proximal operators, but the linear transformation inside $x \mapsto \mathsf{F}(Kx)$ might make the function difficult to optimize directly—for example, we might have $\mathsf{F}(w) = \|w\|_1$ with $K = \nabla_{\mathrm{2d}}$ chosen so that $\mathsf{F}(Kx) = \|x\|_{\mathsf{TV}}$, the two-dimensional total variation norm of $x$. For this type of setting, the Alternating Direction Method of Multipliers (ADMM) reframes the problem as

$$(x^\star, u^\star) = \arg\min_{x,u}\{\mathsf{F}(u) + \mathsf{G}(x) : Kx = v\}$$

and solves for $(x^\star, u^\star)$ by working with the augmented Lagrangian

$$\min_{x,u}\max_{\Delta} L(x,u,\Delta) \text{ where } L(x,u,\Delta) = \left\{\mathsf{F}(u) + \mathsf{G}(x) + \langle \Delta, Kx - u\rangle + \frac{\rho}{2}\|Kx - u\|_2^2\right\}$$

(here $\Delta \in \mathbb{R}^m$ is the dual variable whose role is to enforce the constraint $Kx = u$). The steps of the algorithm are, for each $t \geq 1$,

$$\begin{cases} x_t = \arg\min_x \{L(x, u_{t-1}, \Delta_{t-1})\}, \\ u_t = \arg\min_v \{L(x_t, u, \Delta_{t-1})\}, \\ \Delta_t = \Delta_{t-1} + \rho(Kx_t - u_t). \end{cases}$$

Examining the update step for for $v$, we see that this step entails a single use of the proximal map for $\mathsf{F}$. However, the $x$ update step is more complicated due to the linear operator $K$; this step cannot be solved with one use of the proximal map for $\mathsf{G}$, except in the special case that $K^\top K = \mathbf{I}$. To resolve this, we can add additional curvature to the $x$ update step:

$$x_t = \arg\min_x \left\{L(x, u_{t-1}, \Delta_{t-1}) + \frac{1}{2}(x - x_{t-1})^\top(\lambda\mathbf{I} - \rho K^\top K)(x - x_{t-1})\right\},$$

where $\lambda \geq \rho\|K\|^2$. In this setting, the $x$ update step now becomes solveable with the proximal map for $\mathsf{G}$. In fact, this preconditioned form of the ADMM algorithm is equivalent to the CP algorithm (3); for details of the equivalence, see [7].

In addition to the ADMM and CP algorithms, and their variants, we mention one other option here. If $\mathsf{F}$ is differentiable, then proximal gradient descent offers a simple procedure, alternating between taking a gradient descent step on the term $\mathsf{F}(Kx)$ and a proximal operator step on $\mathsf{G}$. As we discussed in Section 1.1.1, the proximal gradient descent algorithm can be viewed as a simple special case of MOCCA (in that section, the terms were arranged slightly differently, with both differentiable and convex terms all included in $\mathsf{G}$, but the scenarios are equivalent). Of course, this type of method cannot handle scenarios with a nondifferentiable $\mathsf{F}$, which can arise through total variation penalties and in other settings.

### 3.2 The issue of nonconvexity

Next, we turn to the setting where $\mathsf{F}$ and/or $\mathsf{G}$ may be nonconvex.

To begin, we consider the following interesting scenario, introduced in Section 1: suppose that $\mathsf{G}$ is itself convex, with trivial approximations $\mathsf{G}_z = \mathsf{G}$, and that $x \mapsto \mathsf{F}(Kx)$ is strictly convex even though $\mathsf{F}$ is nonconvex. Since $x \mapsto \mathsf{F}(Kx) + \mathsf{G}(x)$ is therefore strictly convex, the optimization problem has a unique global minimizer $x^\star \in \mathbb{R}^d$. However, since $\mathsf{F}$ itself is nonconvex, then the strategies for optimization described in Section 3.1 may not be directly applicable for the task of finding $x^\star$.

Here we outline the difficulties faced by the main existing approaches outlined in Section 3.1 for settings where $\mathsf{F}$ and/or $\mathsf{G}$ are nonconvex (including the scenario outlined above), and summarize the most relevant results in the literature.

To see this, first consider the ADMM algorithm. Suppose that $\mathsf{F}$ is nonconvex, and in particular, for some vector $w \in \mathbb{R}^m$, the function $s \mapsto \mathsf{F}(s \cdot w)$ is strongly concave, in particular,

$$\mathsf{F}(s \cdot w) \leq C - c \cdot s^2 \tag{20}$$

for some $C < \infty, c > 0$ and all $s > 0$. Then we see that the update step

$$u_t = \arg\min_v \{L(x_t, u, \Delta_{t-1})\}$$

is not well-defined; the function $u \mapsto L(x_t, u, \Delta_{t-1})$ diverges to $-\infty$ when we set $u = s \cdot w$ and let $s \to \infty$. Therefore, for some types of nonconvex functions, the ADMM algorithm will not be implementable due to this divergence.

In the literature, theoretical guarantees for the performance of ADMM on nonconvex objective functions have been considered under several different settings. Broadly speaking, we can summarize existing work as falling into one of two categories. First, there are settings where the original form of the ADMM updates perform well. For instance, [24] proves convergence results (for a slightly different algorithm) when optimization is over a bounded set, thus avoiding the issue of divergence arising from directions of strong concavity in $\mathsf{F}$ as mentioned above; this paper also proves results guaranteeing optimality of any limit point, if one exists, in the unbounded optimization setting, but does not guarantee that a limit point is reached. In [19], convergence results are proved assuming that one of the two terms (i.e. $\mathsf{F}$ or $\mathsf{G}$) is smooth; in contrast, for the applications considered here, including total variation type penalties incorporated into $\mathsf{F}$ or hard constraints in $\mathsf{G}$, it is critical to allow for nondifferentiable $\mathsf{F}$. The special case of applying ADMM to nonconvex consensus problems is considered by [15], with convergence guarantees again in a bounded setting, in this case assuming that any nonconvex functions must obey a lower bound.

Second, when the original ADMM updates cannot be expected to perform well—if for instance $\mathsf{F}$ has directions of strong concavity—then the ADMM can be modified by adding curvature to each update via a Bregman divergence term, as studied by Wang et al. [39] (with extensions to multi-block ADMM, with more than two terms in the objective function [40]). This work proves convergence guarantees for the algorithm, but requires that the function $\mathsf{F}$ (after converting to our notation) is differentiable and smooth, in contrast to our work where allowing for nondifferentiable $\mathsf{F}$ is critical.

Next, consider the CP algorithm. When $\mathsf{F}$ is nonconvex, it is no longer the case in general that $\mathsf{F}(Kx) = \max_w\{\langle Kx, w\rangle - \mathsf{F}^*(w)\}$. In fact, by definition of conjugate functions, this maximum defines the "conjugate of the conjugate", i.e.

$$\mathsf{F}^{**}(Kx) = \max_w\{\langle Kx, w\rangle - \mathsf{F}^*(w)\} .$$

It is known that any conjugate function must be convex, i.e. $\mathsf{F}^{**}$ is convex. This implies that $\mathsf{F}^{**} \neq \mathsf{F}$. Therefore, the saddle-point problem does not correspond to the original optimization problem: we have

$$\min_x \max_w\{\langle Kx, w\rangle - \mathsf{F}^*(w) + \mathsf{G}(x)\} = \min_x\{\mathsf{F}^{**}(Kx) + \mathsf{G}(x)\} ,$$

which is different from the original optimization problem since $\mathsf{F}^{**} \neq \mathsf{F}$. If $\mathsf{F}^*(w)$ and $\mathsf{F}^{**}(w)$ take finite values on some domain, then the CP algorithm can be expected to converge, but it will converge to the solution of an optimization problem that is different from the one intended due to the issue that $\mathsf{F}^{**} \neq \mathsf{F}$. Problems also arise in a setting where $\mathsf{F}$ exhibits negative curvature as in (20), in which case $\mathsf{F}^*(w) = \infty$ for all $w$, so we do not have a well-defined saddle point problem to begin with.

To our knowledge, no general results exist for the CP algorithm with nonconvex and nondifferentiable $\mathsf{F}$ and/or $\mathsf{G}$. Valkonen [38] considers an interesting related problem, namely, a variant of the CP algorithm, where the objective function is now $\mathsf{F}(\mathcal{K}(x)) + \mathsf{G}(x)$, for convex $\mathsf{F}, \mathsf{G}$ but with a (nonlinear) map $\mathcal{K}(x)$ in place of the previous linear map $Kx$, as the argument to $\mathsf{F}$. In this case, convergence to a stationary point is proved, even when the nonlinearity of $\mathcal{K}(x)$ may make the overall problem nonconvex. Relatedly, Ochs et al. [29] studies the setting where $\mathsf{F}(\mathcal{K}(x))$ is a nonconvex elementwise penalty on the convex transform $\mathcal{K}(x)$ while $\mathsf{G}$ is convex; their approach uses the CP algorithm as a subroutine for solving convex approximations of the objective function, at each step.

Next, we consider the option of proximal gradient descent, in the case that $\mathsf{G}$ has a simple proximal operator. Loh and Wainwright [22] studies penalized likelihood problems,

$$\min_x \{\mathcal{L}(x) + \mathsf{Penalty}(x)\} ,$$

where the likelihood term $\mathcal{L}(x)$ and/or the penalty term $\mathsf{Penalty}(x)$ may exhibit nonconvexity. In many settings that arise in high-dimensional statistics, for instance, the likelihood term $\mathcal{L}(x)$ may be strongly concave in some directions, but will be strongly convex in all "plausible" directions, that is, all directions $x$ that are not prohibited by the penalty term $\mathsf{Penalty}(x)$. For instance, if $\mathsf{Penalty}(x)$ is a sparsity-promoting penalty, with low values only at solutions $x$ with many (near-)zero values, then $\mathcal{L}(x)$ might be strongly convex in all sparse directions. This relates to the notion of *restricted strong convexity*, introduced by Negahban et al. [26], which we discuss in greater detail in Section 4.2.1. Under restricted convexity and smoothness assumptions on the likelihood term, and with bounds on the amount of nonconvexity allowed in the penalty term, Loh and Wainwright [22] proves convergence to a point that is near the global optimum, for a proximal gradient descent method, with some additional details restricting steps to a bounded set to avoid diverging towards directions of strong concavity. Ochs et al. [28]'s iPiano method gives an inertial (i.e. accelerated) proximal gradient descent method for this same setting where the loss is differentiable with the penalty has an easy proximal map. The (accelerated) proximal gradient descent method for nonconvex problems is studied also by [13, 20]. Note that these algorithms are applicable only when the terms in the objective function are all either differentiable (the likelihood) or have an easy-to-compute proximal operator (the penalty), and therefore, cannot be applied to many of the problems that that we have considered.

Finally, Bolte et al. [3] propose an algorithm, PALM, to solve a related problem of the form

$$\min_{x,w} \{\mathsf{F}(w) + \mathsf{G}(x) + \mathsf{H}(w,x)\}, \tag{21}$$

where $h$ is differentiable while $f, g$ each have easy to compute proximal maps (and may be nonconvex); this formulation is related to the problem we study, but PALM cannot be used to solve general problems of the form $\mathsf{F}(Kx) + \mathsf{G}(x)$ where $\mathsf{F}$ is nondifferentiable (since, if we add a variable $w = Kx$, the constraint $w = Kx$ cannot be enforced with any differentiable function $\mathsf{H}(w,x)$ unless we allow modifications such as a relaxation to a penalty on $\|w - Kx\|_2^2$), and therefore again cannot be applied to some of the problems considered here.

## 3.3 Connection between MOCCA and ADMM

In the convex setting, the CP method with parameter $\theta = 1$ is known to be equivalent to a preconditioned ADMM algorithm [7]. Specifically, reformulating the original optimization problem in the ADMM form

$$\min_{x,u} \{\mathsf{F}(u) + \mathsf{G}(x) : u = Kx\} = \min_{x,u} \max_{\Delta} \left\{ \mathsf{F}(u) + \mathsf{G}(x) + \langle \Delta, Kx - u \rangle + \frac{1}{2} \|Kx - u\|_\Sigma^2 \right\} ,$$

the CP iterations given in (3) are equivalent to the following preconditioned ADMM iterations, where we choose a preconditioning matrix $\mathrm{T}^{-1} - K^\top \Sigma K \succeq 0$:

$$\begin{cases} x_{t+1} = \arg\min_x \left\{ \mathsf{G}(x) + \langle K^\top \Delta_t, x\rangle + \frac{1}{2}\|Kx - u_t\|_\Sigma^2 + \frac{1}{2}\|x - x_t\|_{\mathrm{T}^{-1} - K^\top \Sigma K}^2 \right\}, \\ \Delta_{t+1} = \Delta_t + \Sigma(Kx_{t+1} - u_t), \\ u_{t+1} = \arg\min_u \left\{ \mathsf{F}(u) - \langle \Delta_{t+1}, u\rangle + \frac{1}{2}\|Kx_{t+1} - u\|_\Sigma^2 \right\}, \end{cases}$$

(Here we use a slightly nonstandard indexing, writing the variable updates in the order $x, \Delta, u$ for later convenience.) We do not derive the equivalence here (see [7] for details), but remark that the variables $(x_t, w_t)$ at iteration $t$ of the CP algorithm (3) can be recovered by taking $(x_t, \Sigma(Kx_t - u_t) + \Delta_t)$ from the ADMM iterations.

Similarly, in the more general nonconvex setting considered here, we can equivalently formulate the MOCCA method as a combination of the preconditioned ADMM iterations and taking convex expansions to $\mathsf{F}$ and $\mathsf{G}$. We defer the details to Appendix B, but point out that the ADMM formulation gives us a clearer understanding of the choice of the expansion points $(z, v)$ in the MOCCA algorithm. Recalling the simpler form of the MOCCA algorithm given in Algorithm 1 where the expansion points are updated at each iteration (i.e. there is no "inner loop"), the expansion points were defined as

$$z_{t+1} = x_t \text{ and } v_{t+1} = \Sigma^{-1}(w_{t-1} - w_t) + K\bar{x}_t \ .$$

Using our conversion between the CP variables $(x, w)$ and the ADMM variables $(x, \Delta, v)$, we see that

$$\begin{aligned} v_{t+1} &= \Sigma^{-1}(w_{t-1} - w_t) + K\bar{x}_t \\ &= \Sigma^{-1}(w_{t-1} - w_t) + K(2x_t - x_{t-1}) \quad \text{since we choose extrapolation parameter } \theta = 1 \\ &= \Sigma^{-1}\left(\Sigma(Kx_{t-1} - u_{t-1}) + \Delta_{t-1} - \Sigma(Kx_t - u_t) - \Delta_t\right) + K(2x_t - x_{t-1}) \\ &= \Sigma^{-1}\left(\Sigma(Kx_{t-1} - u_{t-1}) + \Delta_{t-1} - \Sigma(Kx_t - u_t) - (\Delta_{t-1} + \Sigma(Kx_t - u_{t-1}))\right) + K(2x_t - x_{t-1}) \\ &= u_t \ , \end{aligned}$$

where the next-to-last step uses the definition of the update step for $\Delta_t$. In other words, after the $t$th step, our estimated minimizers for $\{\mathsf{F}(u) + \mathsf{G}(x) : Kx = u\}$ are given by $u_t$ and $x_t$, and our convex approximations to $\mathsf{F}$ and to $\mathsf{G}$ are consequently taken at the values $v = u_t$ and $z = x_t$.

# 4 Theoretical results

## 4.1 Convergence to a critical point

First, we show that if the algorithm converges, then its limit point is a solution to our original problem.

**Theorem 1.** *Assume that the families of approximations $\mathsf{F}_v$ and $\mathsf{G}_z$ satisfy* (4)*, and furthermore that*

$$\begin{cases} (v, w) \mapsto \nabla(\mathsf{F} - \mathsf{F}_v)(w) \text{ is continuous jointly in } (v, w), \text{ and} \\ (z, x) \mapsto \nabla(\mathsf{G} - \mathsf{G}_z)(x) \text{ is continuous jointly in } (z, x). \end{cases} \tag{22}$$

*Suppose that Algorithm 2 converges to a point, with*

$$\begin{aligned} &x_{t;\ell} \to \widehat{x} \quad \text{where the sequence } (x_{t;\ell}) \text{ is interpreted as } (x_{1;0}, x_{1;1}, \dots, x_{1;L_1}, x_{2;0}, \dots), \\ &w_{t;\ell} \to \widehat{w} \quad \text{where } (w_{t;\ell}) \text{ is interpreted analogously}, \\ &z_t \to \widehat{z}, \\ &v_t \to \widehat{v}. \end{aligned}$$

*Then $\widehat{x}$ is a critical point of the original optimization problem, in the sense that*

$$0 \in K^\top \partial \mathsf{F}_{K\widehat{x}}(K\widehat{x}) + \partial \mathsf{G}_{\widehat{x}}(\widehat{x}) \ .$$

## 4.2 Guarantees of convergence

We now turn to theoretical results proving that the algorithm converges (and proving rates of convergence) under specific assumptions on $\mathsf{F}$ and $\mathsf{G}$. In this section, we only consider the "inner loop" form of MOCCA, given in Algorithm 2. We show that if our inner loop length (i.e. $L_t$) tends to infinity, then we can give a tight bound on the error of the algorithm.

We begin with an assumption on the step size parameters:

**Assumption 1.** *The extrapolation parameter is set at $\theta = 1$, and the diagonal matrices $\Sigma$ and $\mathrm{T}$ are chosen such that*

$$M = \begin{pmatrix} \mathrm{T}^{-1} & -K^\top \\ -K & \Sigma^{-1} \end{pmatrix} \succ 0 \ .$$

Pock and Chambolle [32] introduce this assumption for the (convex) preconditioned Chambolle-Pock algorithm, and give a simple construction for one choice of $\Sigma, \mathrm{T}$ to satisfy this without calculating any matrix norms or other high-cost operations (see (8)).

Next, we turn to the convexity and smoothness assumptions required for our convergence guarantee.

### 4.2.1 Restricted convexity and smoothness

In practice, $\mathsf{F}$ and/or $\mathsf{G}$ may each consist of multiple terms, combining characteristics of the problem such as a likelihood calculation or a penalty or constraint on the underlying signal. To accomodate a range of potential applications, in particular those arising in the regression and imaging applications described in Section 2, we consider a broad setting where our main assumptions involve the interplay between convexity and negative curvature in the functions $\mathsf{F}, \mathsf{G}$.

The notion of restricted strong convexity (RSC), introduced by Negahban et al. [26], has often been used in high-dimensional statistics to express the idea that likelihood functions and optimization problems, which may not have desirable strong convexity properties globally, nonetheless exhibit strong convexity in "directions of interest". For example, in a least-squares regression problem with design matrix $A \in \mathbb{R}^{n\times d}$, with $n \ll d$, the least squares loss function $\mathcal{L}(x) = \frac{1}{2}\left\|y - Ax\right\|_2^2$ is not strongly convex since $A^\top A$ is rank deficient, but can yield good statistical properties if $A^\top A$ is strongly convex in all sparse directions, that is, $x^\top A^\top Ax \geq c \cdot \left\|x\right\|_2^2$ for all sparse (or approximately sparse) vectors $x$. In this case, the loss function is globally convex, but it is the RSC property that ensures high accuracy for sparse regression problems. More recently, Loh and Wainwright [22] proved that the RSC property, along with an analogous restricted smoothness property, can in fact be leveraged even in nonconvex optimization problems, such as the regression-with-errors-in-variables scenario described in Section 2.1. Their work relies on optimizing the variable $x$ within some bounded set, to ensure that the RSC property will push $x$ towards a good (local) minimum rather than allowing $x$ to diverge. For instance, if the loss function has some directions of strong concavity (as is the case for regression-with-errors-in-variables in (14)), then staying within a bounded set is critical. In theory, their work focuses on problems that take the form of minimizing a penalized loss function over a bounded set $\{x : \left\|x\right\|_1 \leq R\}$, where we think of $R$ as a large bound, requiring only a loose bound on the $\ell_1$ norm of the true signal. In practice, if an optimization algorithm is initialized at zero, then it is often the case that the iterations will never leave a bounded region, without imposing any explicit constraint.

In general, results using the RSC condition take the following form: first, the loss function or objective function $\mathcal{L}(x)$ is shown to satisfy a RSC property of the form

$$\langle x - x', \nabla\mathcal{L}(x) - \nabla\mathcal{L}(x')\rangle \geq c\left\|x - x'\right\|_2^2 - \tau^2\left\|x - x'\right\|_{\mathsf{restrict}}^2 ,$$

for some structured norm $\left\|\cdot\right\|_{\mathsf{restrict}}$ (for example, the $\ell_1$ norm). Here $c > 0$ is a constant while $\tau$ is vanishingly small, for instance $c \sim 1$ and $\tau \sim \sqrt{\frac{\log(d)}{n}}$ in many high dimensional regression applications. The solution $\widehat{x}$ is then shown to converge to the true signal $x^\star$ up to an error of size $\tau R$, where $R$ is some bound on the signal complexity, for instance $R \sim \sqrt{k}$ where $k$ is the true sparsity level of a sparse regression problem. In these settings, it is assumed that $\tau R = o(1)$, and that errors of this magnitude are negligible. We will follow this general framework in our convergence guarantee as well. However, since we consider settings where the signals may not have natural sparsity but would instead have a different type of structure (such as total variation sparsity), we replace the $\ell_1$ norm with a general measure of signal complexity, $\left\|\cdot\right\|_{\mathsf{restrict}}$, chosen with respect to the problem at hand (for instance, a total variation norm).

We now specify our assumptions on convexity and smoothness for the functions involved in the optimization, using the restricted strong convexity / restricted smoothness framework from the literature. Roughly speaking, the following assumption requires that the errors of the convex approximations $\mathsf{F} - \mathsf{F}_v$, $\mathsf{G} - \mathsf{G}_z$ are counterbalanced by strong convexity in the composite approximations $\mathsf{F}_v(Kx) + \mathsf{G}_z(x)$. For the term $\mathsf{G}$, we allow for some flexibility by considering restricted strong convexity and restricted smoothness, relative to the structured norm $\left\|x\right\|_{\mathsf{restrict}}$.

**Assumption 2.** *The approximations* $\mathsf{F}_v$ *and* $\mathsf{G}_z$ *satisfy the conditions* (4)*, with additional assumptions as follows. For the function* $\mathsf{F}$ *and its family of local approximations* $\mathsf{F}_v$*, we assume that* $\mathsf{F}_v$ *is strongly convex, while* $\mathsf{F} - \mathsf{F}_v$ *is smooth: for all* $v, w, u$*,*[3]

$$\begin{cases} \textit{Strong convexity of } \mathsf{F}_v\textit{:}\ \langle u, \partial \mathsf{F}_v(w+u) - \partial \mathsf{F}_v(w)\rangle \geq \|u\|^2_{\Lambda_\mathsf{F}}\,, \\ \textit{Smoothness of } \mathsf{F} - \mathsf{F}_v\textit{:}\ \ |\langle u, \nabla(\mathsf{F} - \mathsf{F}_v)(v+w)\rangle| \leq \frac{1}{2}\left(\|u\|^2_{\Theta_\mathsf{F}} + \|w\|^2_{\Theta_\mathsf{F}}\right), \end{cases}$$

*for some* $\Lambda_\mathsf{F}, \Theta_\mathsf{F} \succeq 0$.[4] *We also assume a gradient condition on* $\mathsf{F}_v$*,*

$$\Big\{\mathsf{F}_v \textit{ satisfies a gradient condition: } \|\partial \mathsf{F}_v(w) - \partial \mathsf{F}_v(w')\|_2 \leq C_{\mathsf{Lip}} + C_{\mathsf{grad}}\,\|w - w'\|_2,$$

*for some* $C_{\mathsf{Lip}}, C_{\mathsf{grad}} < \infty$*. (For example, this is satisfied if* $\mathsf{F}_v$ *can be written as the sum of a Lipschitz function and a smooth function.)*

*For the function* $\mathsf{G}$ *and its family of local approximations* $\mathsf{G}_z$*,we assume that* $\mathsf{G}_z$ *satisfies restricted strong convexity, while* $\mathsf{G} - \mathsf{G}_z$ *satisfies a restricted smoothness assumption: for all* $x, y, z$*,*

$$\begin{cases} \textit{Restricted strong convexity of } \mathsf{G}_z\textit{:}\ \langle y, \partial \mathsf{G}_z(x+y) - \partial \mathsf{G}_z(x)\rangle \geq \|y\|^2_{\Lambda_\mathsf{G}} - \tau^2\,\|y\|^2_{\mathsf{restrict}}\,, \\ \textit{Restricted smoothness of } \mathsf{G} - \mathsf{G}_z\textit{:}\ \ |\langle y, \nabla(\mathsf{G} - \mathsf{G}_z)(z+x)\rangle| \leq \frac{1}{2}\left(\|x\|^2_{\Theta_\mathsf{G}} + \|y\|^2_{\Theta_\mathsf{G}}\right) + \frac{\tau^2}{2}\left(\|x\|^2_{\mathsf{restrict}} + \|y\|^2_{\mathsf{restrict}}\right), \end{cases}$$

*for some* $\Lambda_\mathsf{G}, \Theta_\mathsf{G} \succeq 0$ *and* $\tau < \infty$*.*

*Finally, the total convexity in the local approximations* $\mathsf{F}_v$ *and* $\mathsf{G}_z$ *must (approximately) outweigh the total curvature of the differences* $\mathsf{F} - \mathsf{F}_v$ *and* $\mathsf{G} - \mathsf{G}_z$*. Specifically, for all* $x \in \mathbb{R}^d$*, we require*

$$x^\top(K^\top \Lambda_\mathsf{F} K + \Lambda_\mathsf{G})x \geq x^\top(K^\top \Theta_\mathsf{F} K + \Theta_\mathsf{G})x + C_{\mathsf{cvx}}\,\|x\|_2^2 - \tau^2\,\|x\|^2_{\mathsf{restrict}}\ ,$$

*for some* $C_{\mathsf{cvx}} > 0$ *and* $\tau < \infty$*.*

In general, greater convexity (i.e. $\Lambda_\mathsf{F}, \Lambda_\mathsf{G}$ as strongly positive definite as possible) and tighter bounds on smoothness (i.e. $\Theta_\mathsf{F}, \Theta_\mathsf{G}$ as small as possible) allow for a better (i.e. larger) constant $C_{\mathsf{cvx}}$ and, therefore, faster convergence of the algorithm. The value of $\tau$ is typically of a very small order in many problems arising in high-dimensional statistics, as discussed above.

It is critical to note that this assumption does *not* require either $\mathsf{F}_v$ or $\mathsf{G}_z$ to be strictly convex—if the matrices $\Lambda_\mathsf{F}$ or $\Lambda_\mathsf{G}$ are not full rank, then strict convexity has not been assumed. Instead, $\mathsf{F}_v$ is strongly convex in any direction of $\mathbb{R}^m$ contained in the column span of $\Lambda_\mathsf{F}$, and similarly for $\mathsf{G}_z$ and $\Lambda_\mathsf{G}$ in $\mathbb{R}^d$. Our assumption essentially requires that the combination of these directions leads to overall (approximate) convexity, after accounting for concavity that might be introduced by the errors $\mathsf{F} - \mathsf{F}_v$ and $\mathsf{G} - \mathsf{G}_z$.

For simplicity in the statements and proofs of our results, we group the norms of all matrices from Assumptions 1 and 2 into a single constant:

$$C_{\mathsf{matrix}} = \max\left\{\|\Lambda_\mathsf{F}\|, \|\Theta_\mathsf{F}\|, \|\Lambda_\mathsf{G}\|, \|\Theta_\mathsf{G}\|, \|M\|, \left\|M^{-1}\right\|\right\}\ .$$

Throughout, we will treat $C_{\mathsf{matrix}}$, $C_{\mathsf{cvx}}$, $C_{\mathsf{Lip}}$, and $C_{\mathsf{grad}}$ as fixed finite positive constants, and dependence on these values will not be given explicitly except in the proofs. On the other hand, the role of the restricted convexity/smoothness parameter $\tau$ will be shown explicitly.

### 4.2.2 Convergence guarantee

Choose any point $x^\star \in \mathbb{R}^d$ with $\|x^\star\|_{\mathsf{restrict}} \leq R$, which is a critical point for the optimization problem

$$\min_{\|x\|_{\mathsf{restrict}} \leq R} \{\mathsf{F}(Kx) + \mathsf{G}(x)\}\ .$$

[3] We implicitly restrict all variables to the domain of the appropriate function, throughout.

[4] For a function $\mathsf{F}$ with a multivalued subdifferential, this notation is taken to mean that the statement must hold for any elements of the subdifferentials, throughout the paper.

For convenience, we will now absorb the constraint $\|x\|_{\mathsf{restrict}} \leq R$ into the functions themselves, by replacing $\mathsf{G}$ with the function

$$x \mapsto \begin{cases} \mathsf{G}(x), & \text{if } \|x\|_{\mathsf{restrict}} \leq R, \\ +\infty, & \text{if } \|x\|_{\mathsf{restrict}} > R, \end{cases}$$

and replacing $\mathsf{G}_z$ (for each $z$) with the function

$$x \mapsto \begin{cases} \mathsf{G}_z(x), & \text{if } \|x\|_{\mathsf{restrict}} \leq R, \\ +\infty, & \text{if } \|x\|_{\mathsf{restrict}} > R, \end{cases}$$

In practice, as mentioned before, we typically do not need to explicitly incorporate this constraint into the optimization algorithm, as we will generally only see updates that all lie within a bounded region. However, in our statements and proofs of theoretical results from this point onward, we will assume that $\mathsf{G}, \mathsf{G}_z$ restrict the domain of the variable $x$, that is,

$$\mathsf{G}(x) = \mathsf{G}_z(x) = +\infty \text{ whenever } \|x\|_{\mathsf{restrict}} > R \; . \tag{23}$$

We will also assume that $\tau R$ is bounded by a constant without further comment; since our results give convergence guarantees up to the accuracy level $\tau R$, the results are meaningful only if $\tau R$ is small.

We now state our convergence guarantee for the stable form of the MOCCA algorithm, given in Algorithm 2:

**Theorem 2.** *Assume that Assumptions 1 and 2, and that* $\mathsf{G}, \mathsf{G}_z$ *satisfy* (23). *Then there exists constants* $C_{\mathsf{converge}}, L_{\min} < \infty$ *and* $\delta > 0$*, such that if* $\min_{t\geq 1} L_t \geq L_{\min}$*, then for all* $t \geq 1$*, the iterations of Algorithm 2*

$$\|x_t - x^\star\|_2 \leq C_{\mathsf{converge}} \left( \frac{1}{\sqrt{L'_t}} + \tau R \right)$$

*where*

$$L'_t = \min\left\{L_t, (1+\delta)L_{t-1}, \dots, (1+\delta)^{t-1}L_1\right\} \; .$$

As an example, we can set $L_t \sim (1+\delta)^t$. Then after the $t$th inner loop, $\|x_t - x^\star\|_2 \sim (1+\delta)^{-t/2} + \tau R$, and the total number of iterations taken is $L_1 + \dots + L_t \sim (1+\delta)^t$. In other words, the error $\|x - x^\star\|_2$ scales as $\frac{1}{\sqrt{T}} + \tau R$ where $T$ is the total number of update steps, i.e. error is inversely proportional to the square root of computational cost, up to the accuracy level $\tau R$.

*Remark* 1. If we were to assume additionally that $\mathsf{F}$ is differentiable and smooth, the result would improve dramatically: we would obtain error decaying exponentially in the number of update steps (up to the accuracy level $\tau R$). The resulting convergence guarantee would then be comparable to the results obtained in Loh and Wainwright [22, Theorem 3], which show a result with error at time $t$ scaling as $c^t + \tau R$ (for a constant $c < 1$). However, convergence in this setting is of limited interest for the applications we have in mind, since total variation penalties, and many other natural penalties or losses falling into the $\mathsf{F}$ term of the composite objective functions, are not differentiable.

### 4.2.3 Convexity and smoothness assumptions: an example

To illustrate the many different matrices and constants appearing in Assumption 2 with a concrete example, we return to the problem studied in Section 2.1, where a least squares regression with errors in variables is combined with a total variation (or generalized $\ell_1$) penalty. Recalling this setting, we seek to minimize $\mathcal{L}(x) + \nu \cdot \|Kx\|_1$, and we set $\mathsf{F}_v(w) = \mathsf{F}(w) = \nu \cdot \|w\|_1$,

$$\mathsf{G}(x) = \begin{cases} \mathcal{L}(x), & \text{if } \|x\|_1 \leq R, \\ +\infty, & \text{if } \|x\|_{\mathsf{restrict}} > R, \end{cases}$$

and

$$\mathsf{G}_z(x) = \begin{cases} \mathcal{L}(x) + \frac{\sigma_A^2}{2}\|x - z\|_2^2, & \text{if } \|x\|_{\mathsf{restrict}} \leq R, \\ +\infty, & \text{if } \|x\|_1 > R, \end{cases}$$

where

$$\mathcal{L}(x) = \frac{1}{2}x^\top \left(\frac{Z^\top Z}{n} - \sigma_A^2 \mathbf{I}_d\right) x - x^\top \left(\frac{Z^\top b}{n}\right).$$

Under the Gaussian noise model, the noisy design matrix given by entries $Z_{ij} = A_{ij} + \mathsf{Normal}(0, \sigma_A^2)$ and the response is given by $b = A \cdot x_{\mathsf{true}} + \mathsf{Normal}(0, \sigma^2 \mathbf{I})$. The MOCCA update steps for this problem are given in (16).

In this setting, Assumption 2 is satisfied with the following parameters. First, since $\mathsf{F}$ is convex but not strongly convex, we set $\Lambda_{\mathsf{F}} = \mathbf{0}$; we can also set $\Theta_{\mathsf{F}} = \mathbf{0}$ as $\mathsf{F}_v = \mathsf{F}$ for any expansion point $v$. $\mathsf{F}$ is $\nu$-Lipschitz so we can take $C_{\mathsf{Lip}} = \nu, C_{\mathsf{grad}} = 0$. Next, for $\mathsf{G}$, we see that

$$\mathsf{G}_z(x) = \frac{1}{2}x^\top \left(\frac{Z^\top Z}{n}\right) x - x^\top \left(\frac{Z^\top b}{n} + \sigma_A^2 z\right)$$

(on the domain $\|x\|_{\mathsf{restrict}} \le R$), and so we can set $\Lambda_{\mathsf{G}} = \frac{Z^\top Z}{n}$. To check the smoothness condition, we have

$$\langle y, \nabla(\mathsf{G} - \mathsf{G}_z)(z + x)\rangle = \langle y, \sigma_A^2 x\rangle \le \sigma_A^2 \cdot \frac{1}{2}(\|x\|_2^2 + \|y\|_2^2),$$

and so we can take $\Theta_{\mathsf{G}} = \sigma_A^2 \mathbf{I}_d$.

Finally, for the "total convexity" condition of Assumption 2, we need to check that $K^\top(\Lambda_{\mathsf{F}} - \Theta_{\mathsf{F}})K + (\Lambda_{\mathsf{G}} - \Theta_{\mathsf{G}}) = \frac{Z^\top Z}{n} - \sigma_A^2 \mathbf{I}_d$ satisfies restricted strong convexity. In Loh and Wainwright [21, Corollary 1], it is shown that if the rows of the (original) design matrix $A$ are drawn i.i.d. from a subgaussian distribution with covariance $\Sigma_A$ then, assuming that the sample size $n$ satisfies $n \gg \log(d)$, the matrix $\frac{Z^\top Z}{n} - \sigma_A^2 \mathbf{I}_d$ (which is an unbiased estimate of the desired term $\frac{A^\top A}{n}$ using the unknown original design matrix $A$) satisfies

$$x^\top \left(\frac{Z^\top Z}{n} - \sigma_A^2 \mathbf{I}_d\right) x \ge \frac{1}{2}\lambda_{\min}(\Sigma_A) \cdot \|x\|_2^2 - (\text{constant}) \cdot \frac{\log(d)}{n} \cdot \|x\|_{\mathsf{restrict}}^2 \text{ for all } x \in \mathbb{R}^d \tag{24}$$

with high probability, when we choose $\|x\|_{\mathsf{restrict}} = \|x\|_1$; similar results will hold for other structured choices of $\|\cdot\|_{\mathsf{restrict}}$ such as total variation norm or a generalized $\ell_1$ norm. Therefore we can set $C_{\mathsf{cvx}} = \frac{1}{2}\lambda_{\min}(\Sigma_A)$ and $\tau \sim \sqrt{\frac{\log(d)}{n}}$ to obtain the desired condition in Assumption 2. Note that the guarantees of Theorem 2 give a meaningful convergence result even if we choose a fairly large radius $R$.

# 5 Experiments

We now implement the MOCCA algorithm to examine its performance in practice. Throughout this section, we work with the simpler formulation of MOCCA, given in Algorithm 1, with no "inner loop". All computations were performed in MATLAB [25].[5]

For all simulations, we choose not to place a bound on $\|x\|_{\mathsf{restrict}}$, which technically is required by our convergence guarantees and those of the related results in [22] (which we compare to, in Simulation 2). Empirically we observe good convergence without imposing such a bound, but can easily add such a bound if desired.

We consider two examples: Simulation 1 studies nonconvex total variation regularization with a least squares loss (as described in Section 2.3), and Simulation 2 considers convex total variation regularization with a nonconvex loss arising from regression with errors in variables (as described in Section 2.1). While other algorithms which are developed specifically for these problems are available—for example, denoising with total variation penalties is studied by e.g. [6, 42, 43] and could be combined with a proximal gradient method for Simulation 2—here our purpose is simply to illustrate applications of MOCCA to several concrete examples in order to demonstrate its flexibility for a broad range of problems. Specific problems will often have specialized algorithms which would far outperform our general-purpose method; however, slight modifications to the optimization problem (for example, replacing total variation regularization with a more general penalty $\|Kx\|_1$ for a generic dense matrix $K$, or with isotropic total variation) will often mean that specialized algorithms can no longer be applied while MOCCA can adapt easily to accomodate these changes.

[5]Code for fully reproducing these simulations is available at `http://www.stat.uchicago.edu/~rina/mocca.html`.

## 5.1 Simulation 1: nonconvex total variation penalty

In the first simulation, we study the nonconvex total variation penalty considered in Section 2.3, using a two-dimensional spatial structure. We generate data as follows: first, we define the true signal $x_{\mathsf{true}} \in \mathbb{R}^d$ with dimension $d = 625$, obtained by vectorizing the two-dimensional locally constant array

$$\begin{pmatrix} \mathbf{1}_{5\times 5} & \mathbf{0}_{5\times 15} & \mathbf{0}_{5\times 5} \\ \mathbf{0}_{15\times 5} & \mathbf{1}_{15\times 15} & \mathbf{0}_{15\times 5} \\ \mathbf{0}_{5\times 5} & \mathbf{0}_{5\times 15} & \mathbf{1}_{5\times 5} \end{pmatrix} \in \mathbb{R}^{25\times 25}.$$

The two-dimensional total variation of the true signal is very low, because $\nabla_{\mathrm{2d}} x_{\mathsf{true}}$ is sparse. We then take a linear regression model with $n = 200$ observations, with design matrix $A \in \mathbb{R}^{n\times d}$ with $A_{ij} \stackrel{\mathrm{iid}}{\sim} \mathsf{Normal}(0,1)$ and $b \in \mathbb{R}^n$ with entries

$$b_i = (A \cdot x_{\mathsf{true}})_i + \mathsf{Normal}(0,1) \; .$$

We would then like to solve a penalized least-squares problem using the nonconvex total variation penalty introduced in Section 2.3, namely,

$$\widehat{x} = \arg\min_x \{\mathsf{Obj}(x)\} \text{ for } \mathsf{Obj}(x) = \frac{1}{2}\left\|b - Ax\right\|_2^2 + \nu \cdot \mathsf{logTV}_\beta(x) \; , \tag{25}$$

where we choose penalty parameter $\nu = 20$ and nonconvexity parameter $\beta = 3$ (recall that a low value of $\beta$ corresponds to greater nonconvexity), and where the $\mathsf{logTV}_\beta(\cdot)$ penalty is defined with respect to two-dimensional total variation—recall

$$\mathsf{logTV}_\beta(x) = \mathsf{logL1}_\beta(\nabla_{\mathrm{2d}} x) = \sum_i \beta \log\left(1 + |(\nabla_{\mathrm{2d}} x)_i|/\beta\right) \; .$$

Here $\nabla_{\mathrm{2d}} \in \mathbb{R}^{m\times d}$ is the two-dimensional first differences matrix for the vectorized $d_1 \times d_2$ grid, where $d = d_1 \cdot d_2$ is the total dimension of the signal while $m = d_1(d_2 - 1) + d_2(d_1 - 1)$ is the number of first-order differences measured; in our case, we have $d_1 = d_2 = 25$, $d = 625$, and $m = 1200$.

Next, we implement the MOCCA algorithm with the two variants described in Section 2.3: setting $K = \nabla_{\mathrm{2d}}$, we consider the more natural form where the penalty term is contained in $\mathsf{F}$, given by

$$\begin{cases} \mathsf{F}(w) = \nu \cdot \mathsf{logL1}_\beta(w), \text{ with } \mathsf{F}_v(w) = \nu \cdot \left\|w\right\|_1 + \nu \left[h_\beta(v) + \langle w - v, \nabla h_\beta(v)\rangle\right], \\ \mathsf{G}(x) = \mathsf{G}_z(x) = \frac{1}{2}\left\|b - Ax\right\|_2^2, \end{cases} \tag{26}$$

where $h_\beta(w) = \mathsf{logL1}_\beta(w) - \left\|w\right\|_1$ is a differentiable concave function as discussed in Section 2.3. We also consider the less natural form where the penalty term is split across $\mathsf{F}$ and $\mathsf{G}$, given by

$$\begin{cases} \mathsf{F}(w) = \mathsf{F}_v(w) = \nu \cdot \left\|w\right\|_1, \\ \mathsf{G}(x) = \frac{1}{2}\left\|b - Ax\right\|_2^2 + \nu \cdot h_\beta(\nabla_{\mathrm{2d}} x), \\ \qquad \text{with } \mathsf{G}_z(x) = \frac{1}{2}\left\|b - Ax\right\|_2^2 + \nu\left[h_\beta(\nabla_{\mathrm{2d}} z) + \langle \nabla_{\mathrm{2d}}(x - z), \nabla h_\beta(\nabla_{\mathrm{2d}} z)\rangle\right], \end{cases} \tag{27}$$

We will refer to these two versions as MOCCA(natural) and MOCCA(split), respectively. Finally, we choose step size parameters

$$\Sigma = \lambda \cdot \frac{1}{2}\mathbf{I}_m \quad \text{and} \quad \mathrm{T} = \lambda^{-1} \cdot \frac{1}{4}\mathbf{I}_d \; ,$$

which ensures that the positive semidefinite assumption, Assumption 1, will hold (although perhaps not strictly) as in [32]. We test the algorithm across a range of $\lambda$ values, $\lambda \in \{4, 8, 16, 32, 64\}$.

The results are shown in Figure 2, which plots the log value of the objective function $\mathsf{Obj}(x)$ at each iteration, and also plots the log of the change in each iteration,

$$\mathsf{Change}_t = \left\| \begin{pmatrix} x_{t-1} - x_t \\ w_{t-1} - w_t \end{pmatrix} \right\|_2 \; . \tag{28}$$

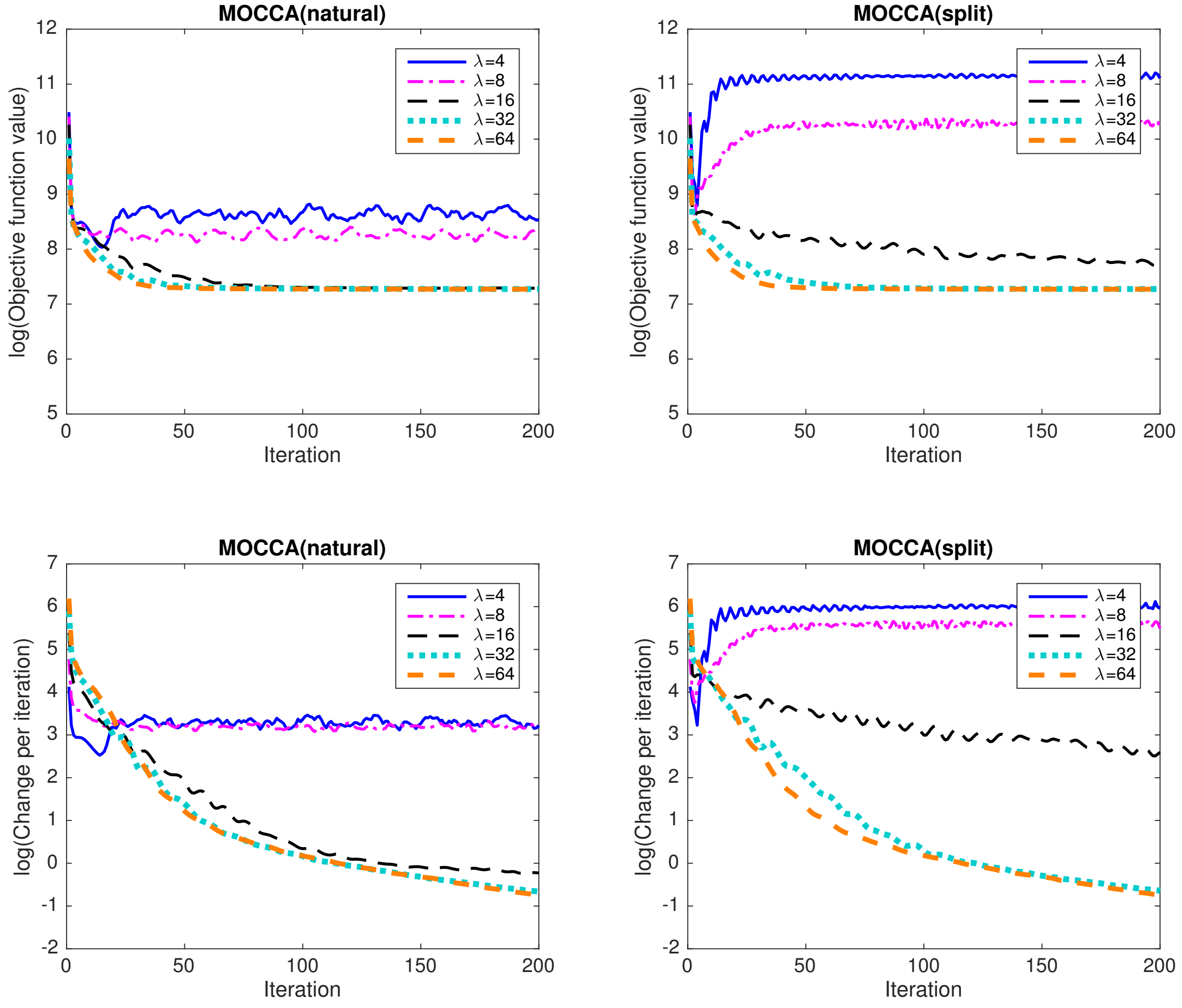

Figure 2: Results from Simulation 1. The top row plots the log of the objective function value defined in (25) against iteration number, while the bottom row plots the log of the change in the $(x, w)$ variables (28) at each iteration, for several step size parameters $\lambda$. (Section 1.1.4 gives a direct correspondence between convergence of the $(x, w)$ variables, and the optimality of the $x$ variable.) The plots show results from the "natural" (left) and "split" (right) versions of the MOCCA algorithm, defined in (26) and (27), respectively.

(Recall from Section 1.1.4 that, if $\mathsf{Change}_t \to 0$, then the optimality gap of the solution tends to zero; that is, the $x$ variable is close to being a critical point for the optimization problem.) Looking first at the results for MOCCA(natural), we see that smaller $\lambda$ values tend to lead to faster convergence at the very early stages, but poorer performance or instability at later stages. (In fact, this suggests the possibility of varying $\lambda$ as we run more iterations, which we leave to future work.)

Turning to MOCCA(split), we see that the performance is worse at all $\lambda$ values as compared with MOCCA(natural); the difference is minor for the largest $\lambda$ values, but the lower $\lambda$ values give far poorer results and far more instability for MOCCA(split) as compared to MOCCA(natural). This highlights the importance of the "mirroring" step in our algorithm, which gives us the flexibility of placing the nonconvex terms into $\mathsf{F}$, i.e. the function which will be optimized via the dual variable. In other scenarios, of course, a different arrangement of the terms may be preferable.

## 5.2 Simulation 2: total variation penalty, with errors in variables

In the second simulation, we treat the errors-in-variables setting discussed in Section 2.1. We generate the signal $x_{\mathsf{true}}$, the design matrix $A$, and the response vector $b$ as in Simulation 1. Next, suppose our measurement of $A$ is itself noisy: define $Z \in \mathbb{R}^{n\times d}$ with

$$Z_{ij} = A_{ij} + \sigma_A \cdot \mathsf{Normal}(0,1) ,$$

where $\sigma_A = 0.2$. Finally, we would like to minimize the objective function

$$\frac{1}{2}x^\top \left(Z^\top Z - n\cdot\sigma_A^2\mathbf{I}_d\right)x - x^\top Z^\top b + \nu\left\|x\right\|_{\mathsf{TV}} , \tag{29}$$

with penalty parameter $\nu = 20$, where again we use two-dimensional total variation, $\|x\|_{\mathsf{TV}} = \|\nabla_{\mathrm{2d}}x\|_1$. (Here we use a different scaling of the likelihood term relative to Section 2.1 for simpler implementation and tuning.) Of course, due to the negative quadratic term, this objective function is strongly concave in some directions and so its global minimum is "at infinity"; within a bounded set, however, the penalty will ensure that the objective function is approximately convex. In practice, initializing our algorithm at $x = 0$ does not lead to any problems, and we converge to a bounded solution that can be viewed as a local minimum within a bounded set, e.g. $\{x : \|x\|_{\mathsf{TV}} \leq R\}$ for some appropriate choice of $R$, as discussed in the context of restricted strong convexity in Section 4.2.1.

A proximal gradient descent method, as proposed by [22] for this type of nonconvex penalized likelihood, in theory would iterate the steps

$$\begin{cases}\tilde{x}_{t+1} = x_t - \frac{1}{\eta}\left(\left(Z^\top Z - n\sigma_A^2\mathbf{I}_d\right)x_t - (Z^\top b)\right) ,\\ x_{t+1} = \arg\min_x\left\{\frac{1}{2}\left\|x - \tilde{x}_{t+1}\right\|_2^2 + \frac{\nu}{\eta}\left\|x\right\|_{\mathsf{TV}}\right\} ,\end{cases}$$

where $\frac{1}{\eta}$ is a step size parameter. However, the second step is a proximal operator for the total variation norm $\|x\|_{\mathsf{TV}} = \|\nabla_{\mathrm{2d}}x\|_1$, which cannot be calculated in closed form. Instead, we could apply the CP algorithm to the convex (sub)problem of this proximity operator with parameters $\Sigma = \lambda\cdot\frac{1}{2}\mathbf{I}_m$ and $\mathrm{T} = \lambda^{-1}\cdot\frac{1}{4}\mathbf{I}_d$, and could terminate this inner "prox loop" after some convergence criterion is reached, e.g. after some fixed number $n_{\mathsf{step}}$ of steps, or once the relative change in $x$ is below $\epsilon_{\mathsf{thresh}}$. We do not show details of the derivation, but the complete procedure iterates these steps:

$$\begin{cases}\text{Gradient step: } \tilde{x}_{t+1} = x_t - \frac{1}{\eta}\left(\left(Z^\top Z - n\sigma_A^2\mathbf{I}_d\right)x_t - (Z^\top b)\right) ,\\ \text{Initialize prox loop: } x'_{t+1;0} = x_t, u'_{t+1;0} = u_t.\\ \text{Run prox loop: for } \ell = 1,2,\dots ,\\ \begin{cases} x'_{t+1;\ell} = (1+\frac{1}{4\lambda})^{-1}\left(x'_{t+1;\ell-1} + \frac{1}{4\lambda}\tilde{x}_{t+1} - \frac{1}{4\lambda}\nabla_{\mathrm{2d}}^\top u'_{t+1;\ell-1}\right) ,\\ u'_{t+1;\ell} = \mathsf{Truncate}_{[-\nu/\eta,\nu/\eta]}\left(u'_{t+1;\ell-1} + \frac{\lambda}{2}\nabla_{\mathrm{2d}}\bar{x}'_{t+1;\ell}\right) \text{ where } \bar{x}'_{t+1;\ell} = x'_{t+1;\ell} + \theta(x'_{t+1;\ell} - x'_{t+1;\ell-1}),\\ \text{until a convergence criterion is reached (i.e. } \ell = n_{\mathsf{step}} \text{ or } \left\|x'_{t+1} - x'_t\right\|_2 / \left\|x'_t\right\|_2 \leq \epsilon_{\mathsf{thresh}}).\end{cases}\\ \text{Gather results from prox loop: } x_{t+1} = x'_{t+1;\ell}, u_{t+1} = u'_{t+1;\ell}.\end{cases} \tag{30}$$

We will refer to this method as the Approximate Proximal Gradient Descent (APGD) algorithm, where "approximate" describes the fact that the proximal operator step is only ever solved approximately via a finite number of steps in the inner loop.

In fact, if we examine this algorithm carefully, we can find that by taking a single step of the inner "prox loop" (that is, setting $n_{\mathsf{step}} = 1$), we arrive back at the steps of the MOCCA algorithm. Specifically, as in the implementation (15) in Section 2.1, we choose $K = \nabla_{\mathrm{2d}}$, $\mathsf{F}(w) = \mathsf{F}_v(w) = \nu\left\|w\right\|_1$, and

$$\mathsf{G}(x) = \frac{1}{2}x^\top\left(Z^\top Z - n\sigma_A^2\mathbf{I}_d\right)x - x^\top(Z^\top b)$$

with local approximations given by linear expansions,

$$\mathsf{G}_z(x) = \mathsf{G}(z) + \langle x - z, \nabla \mathsf{G}(x) \rangle = \langle x, (Z^\top Z - n\sigma_A^2 \mathbf{I}_d) z - Z^\top b \rangle + \text{(terms constant with respect to } x\text{)}.$$

The update steps of the MOCCA algorithm are then given by

$$\begin{cases} x_{t+1} = \arg\min_x \left\{ \langle \nabla_{2\mathrm{d}} x, w_t \rangle + \mathsf{G}_{z_t}(x) + \frac{1}{2} \|x - x_t\|^2_{\mathrm{T}^{-1}} \right\}, \\ w_{t+1} = \arg\min_y \left\{ -\langle \nabla_{2\mathrm{d}} \bar{x}_{t+1}, w \rangle + \mathsf{F}^*(w) + \frac{1}{2} \|w - w_t\|^2_{\Sigma^{-1}} \right\}, \\ z_{t+1} = x_{t+1} , \end{cases}$$

which we can simplify to

$$\begin{cases} x_{t+1} = x_t - \mathrm{T}\left( \nabla_{2\mathrm{d}}^\top w_t + \left(Z^\top Z - n\sigma_A^2 \mathbf{I}_d\right) x_t - (Z^\top b) \right), \\ w_{t+1} = \mathsf{Truncate}_{[-\nu,\nu]}\left( w_t + \Sigma \nabla_{2\mathrm{d}} \bar{x}_{t+1} \right). \end{cases}$$

If we choose

$$\Sigma = \frac{\lambda\eta}{2} \cdot \mathbf{I}_{d-1} \quad \text{and} \quad \mathrm{T} = \frac{1}{(4\lambda+1)\eta} \cdot \mathbf{I}_d ,$$

it can be shown that this is equivalent to the proximal gradient algorithm (30) with a single inner loop step, i.e. with $n_{\mathsf{step}} = 1$ (specifically, the iterates $x_t$ stay the same, while the other variables are related as $w_t = \eta \cdot u_t$).

Now we compare the performance of the approximate proximal gradient descent (APGD) algorithm, with various stopping criteria for the inner "prox loop", against the performance of the MOCCA algorithm, which we can view as the APGD algorithm taking exactly one step in each inner "prox loop". For simplicity, we consider only a few values for the step size parameters, setting $\eta = \lambda = 100$ or $\eta = \lambda = 200$. As for Simulation 1, we will see that higher values for these parameters gives more stability at the cost of slower convergence.

We consider stopping rules for the inner loop as follows: either we run the inner loop for a fixed number of steps, $n_{\mathsf{step}} \in \{1, 5\}$ (with $n_{\mathsf{step}} = 1$ yielding MOCCA), or we use a convergence criterion $\epsilon_{\mathsf{thresh}} \in \{0.1, 0.05, 0.01\}$. Figure 3 shows the results; for the figure on the left, we see that running the inner loop longer does help to make our solutions more accurate (i.e. the objective function is lower) over the range of iterations. However, each iteration has greater computational cost when we increase the time spent running the inner loop. Since we would like to see the performance as a function of computational cost, the plot on the right-hand side of Figure 3 shows the same results plotted against the true number of iterations, i.e. where we count each pass through the inner loop of (30) rather than counting only passes through the outer loop of (30). In this setting, we see that in fact the various versions of the algorithms perform nearly identically—in other words, a one-step approximation to the proximal map performs just as well as a more conservative inner loop that is run for longer—with the exception of setting $\epsilon_{\mathsf{thresh}} = 0.01$ and $\eta = \lambda = 100$, in which case the algorithm diverges immediately. It is interesting to note that this small choice for $\epsilon_{\mathsf{thresh}}$ is closest in spirit to proximal gradient descent (that is, the proximal step is the most accurate), although of course there may be some effect of tuning parameters. The choice $\epsilon_{\mathsf{thresh}}$ does achieve convergence with the more conservative choice $\eta = \lambda = 200$, but convergence is noticeably slower in this case. Thus, we can conclude that when the proximal step does not have a closed form solution, it may be better to use a coarse approximation (which, implicitly, is the strategy taken by MOCCA for this problem) rather than aiming for near-convergence in the proximal step for each iteration.

# 6 Proofs

## 6.1 Critical points (Theorem 1)

For this proof we will use two facts: for a continuous convex function $h$,

$$\text{If } a_t \in \partial h(b_t) \text{ and } a_t \to a,\, b_t \to b \text{ then } a \in \partial h(b), \tag{31}$$

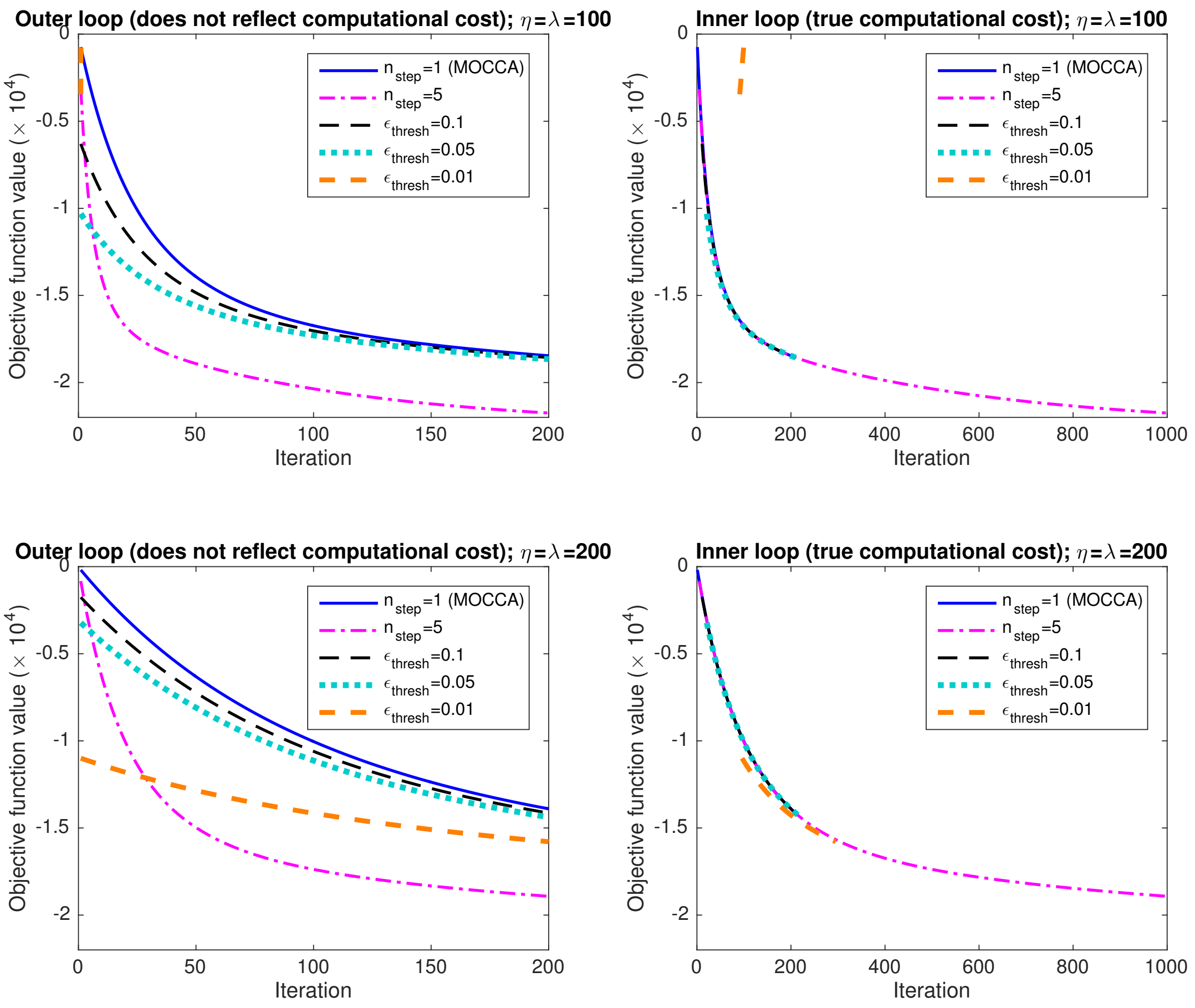

Figure 3: Results from Simulation 2, plotting the value of the objective function (29) against iteration number, counting either iterations of the "outer loop" (left) or of the "inner loop" (right), for various stopping rules for the inner loop in (30). Parameters are set as $\eta = \lambda = 100$ (top) or $\eta = \lambda = 200$ (bottom). Counting the number of passes through the inner loop is an accurate reflection of the true computational cost of (30), and so the right-hand plots give the correct interpretation of the results, where the various versions of the algorithms perform nearly identically in both settings (the lines are indistinguishable) except for the setting $\epsilon_{\mathsf{thresh}} = 0.01$ which diverges for the setting $\eta = \lambda = 100$. Setting $n_{\mathsf{step}} = 1$ yields the MOCCA algorithm as discussed in the text.

and

$$a \in \partial h(b) \text{ if and only if } b \in \partial h^*(a). \tag{32}$$

First, by definition of $w_{t+1} = w_{t;L_t}$,

$$\partial \mathsf{F}^*_{v_t}(w_{t+1;1}) \ni \Sigma^{-1}(w_{t+1;0} - w_{t+1;1}) + K\bar{x}_{t+1;1} \ ,$$

and therefore by (32),

$$w_{t+1;1} \in \partial \mathsf{F}_{v_t}(\Sigma^{-1}(w_{t+1;0} - w_{t+1;1}) + K\bar{x}_{t+1;1}).$$

We can rewrite this as

$$w_{t+1;1} \in \underbrace{\partial \mathsf{F}_{\widehat{v}}(\Sigma^{-1}(w_{t+1;0} - w_{t+1;1}) + K\bar{x}_{t+1;1})}_{\text{(Term 1)}} + \underbrace{\nabla(\mathsf{F}_{v_t} - \mathsf{F}_{\widehat{v}})(\Sigma^{-1}(w_{t+1;0} - w_{t+1;1}) + K\bar{x}_{t+1;1})}_{\text{(Term 2)}}. \tag{33}$$

Next, since the solution converges, we see that

$$w_{t+1;1} \to \widehat{w}, \ v_t \to \widehat{v}, \ \Sigma^{-1}(w_{t+1;0} - w_{t+1;1}) + K\bar{x}_{t+1;1} \to K\widehat{x}.$$

Our assumption (22) implies that

$$(v,w) \mapsto \nabla(\mathsf{F}_v - \mathsf{F}_{\widehat{v}})(w) = -\nabla(\mathsf{F} - \mathsf{F}_v)(w) + \nabla(\mathsf{F} - \mathsf{F}_{\widehat{v}})(w)$$

is jointly continuous in $(v,w)$, and so

$$\text{(Term 2)} = \nabla(\mathsf{F}_{v_t} - \mathsf{F}_{\widehat{v}})(\Sigma^{-1}(w_{t+1;0} - w_{t+1;1}) + K\bar{x}_{t+1;1}) \to \nabla(\mathsf{F}_{\widehat{v}} - \mathsf{F}_{\widehat{v}})(K\widehat{x}) = 0.$$

Therefore, applying the property (31) to the expression in (33), we see that

$$\widehat{w} \in \partial \mathsf{F}_{\widehat{v}}(K\widehat{x}). \tag{34}$$

Next, for each $t$, by definition,

$$\begin{aligned} v_{t+1} = \frac{1}{L_{t+1}} \sum_{\ell=1}^{L_{t+1}} \left(\Sigma^{-1}(w_{t+1;\ell-1} - w_{t+1;\ell}) + K\bar{x}_{t+1;\ell}\right) \\ = \frac{1}{L_{t+1}} \left(\Sigma^{-1}(w_{t+1;0} - w_{t+1;L_{t+1}}) + K(x_{t+1;L_{t+1}} - x_{t+1;0})\right) + Kx_{t+1} \ . \end{aligned}$$

Taking limits on each side, we see that $\widehat{v} = K\widehat{x}$. Returning to (34) above this proves that

$$\widehat{w} \in \partial \mathsf{F}_{K\widehat{x}}(K\widehat{x}) \ .$$

Next, by definition of the algorithm,

$$x_{t+1;1} = (\mathbf{I} + \mathrm{T}\partial \mathsf{G}_{z_t})^{-1}(x_{t+1;0} - \mathrm{T}K^\top w_{t+1;0}) \ ,$$

and after rearranging terms,

$$\partial \mathsf{G}_{z_t}(x_{t+1;1}) = \mathrm{T}^{-1}(x_{t+1;0} - x_{t+1;1}) - K^\top w_{t+1;0} \ .$$

Taking limits on each side as $t \to \infty$, and applying (31),

$$\partial \mathsf{G}_{\widehat{z}}(\widehat{x}) \ni -K^\top \widehat{w} \ .$$

And, we know that $z_{t+1} = x_{t+1}$ for each $t$, therefore $\widehat{z} = \widehat{x}$, and so

$$\partial \mathsf{G}_{\widehat{x}}(\widehat{x}) \ni -K^\top \widehat{w} \ .$$

Combining the work above, then,

$$0 = K^\top \widehat{w} - K^\top \widehat{w} \in K^\top \partial \mathsf{F}_{K\widehat{x}}(K\widehat{x}) + \partial \mathsf{G}_{\widehat{x}}(\widehat{x}) \ ,$$

as desired.

## 6.2 Convergence guarantee (Theorem 2)

### 6.2.1 Notation

Fixing any expansion points $(z,v) \in \mathbb{R}^d \times \mathbb{R}^m$, we define a primal-dual update step:

$$(x', w') = \mathsf{Step}_{z,v}(x,w) \ ,$$

given by

$$\begin{cases} x' \leftarrow (\mathbf{I} + \mathrm{T}\partial \mathsf{G}_z)^{-1}(x - \mathrm{T}^{-1}K^\top w) \ , \\ w' \leftarrow (\mathbf{I} + \Sigma \partial \mathsf{F}^*_v)^{-1}(w + \Sigma K(2x' - x)) \ . \end{cases}$$

This is one step of the CP algorithm applied to the convex objective function

$$\min_x \{\mathsf{F}_v(Kx) + \mathsf{G}_z(x)\}$$

with extrapolation parameter $\theta = 1$.

Next, defining $x^\star$ to be any critical point of the original problem as before, let $w^\star \in \partial \mathsf{F}_{Kx^\star}(Kx^\star)$ be an element of the subdifferential such that

$$0 \in K^\top w^\star + \partial \mathsf{G}_{x^\star}(x^\star) \ .$$

Finally, for any expansion points $(z, v)$, define the (not necessarily unique) solution for the convex optimization problem when the we use approximations $\mathsf{F}_v, \mathsf{G}_z$:

$$x^\star_{z,v} = \arg\min_x \{\mathsf{F}_v(Kx) + \mathsf{G}_z(x)\} \ ,$$

and let $w^\star_{z,v} \in \partial \mathsf{F}_v(Kx^\star_{z,v})$ be an element of the subdifferential such that

$$0 \in K^\top w^\star + \partial \mathsf{G}_z(x^\star_{z,v}) \ .$$

### 6.2.2 Lemmas

The proof of Theorem 2 can be split into several key results. First we state these lemmas and explain their role, then we will formally prove the theorem. The lemmas are proved in Appendix A.

The first lemma shows that, if we use expansion points $(z, v)$ close to the true solution, i.e. $(z, v) \approx (x^\star, Kx^\star)$, then the minimizer $x^\star_{z,v}$ for the convex approximation will be close to $x^\star$.

**Lemma 1.** *Suppose that Assumptions 1 and 2 hold. Define*

$$\Theta = \begin{pmatrix} \Theta_{\mathsf{G}} & 0 \\ 0 & \Theta_{\mathsf{F}} \end{pmatrix} + \frac{C_{\mathsf{cvx}}}{(C_{\mathsf{matrix}})^2}\mathbf{I} \succ 0 \ .$$

*Then there exist constants $C_{\mathsf{contr}} > 0, C_{\mathsf{excess}} < \infty$, which depend only on $C_{\mathsf{matrix}}, C_{\mathsf{cvx}}, C_{\mathsf{Lip}}, C_{\mathsf{grad}}$, such that for any $(z, v)$ with $\|z\|_{\mathsf{restrict}} \leq R$,*

$$\left\|\begin{pmatrix} x^\star_{z,v} - x^\star \\ Kx^\star_{z,v} - Kx^\star \end{pmatrix}\right\|_\Theta \leq (1 - C_{\mathsf{contr}}) \left\|\begin{pmatrix} z - x^\star \\ v - Kx^\star \end{pmatrix}\right\|_\Theta + C_{\mathsf{excess}} \cdot \tau R$$

*and*

$$\left\|\begin{pmatrix} x^\star_{z,v} - x^\star \\ w^\star_{z,v} - w^\star \end{pmatrix}\right\|_2 \leq C_{\mathsf{Lip}} + C_{\mathsf{excess}} \left( \left\|\begin{pmatrix} z - x^\star \\ v - Kx^\star \end{pmatrix}\right\|_\Theta + \tau R \right) \ .$$

The second lemma shows that, after running an "inner loop", the $(x, w)$ variables are nearly optimal for the current convex approximation, and the next expansion points are also near this optimum.

**Lemma 2.** *Suppose that Assumptions 1 and 2 hold. For any $L \geq 1$, and any $x^{(0)}, w^{(0)}, z, u$ with $\|x^{(0)}\|_{\mathsf{restrict}} \leq R$, suppose we iterate $\mathsf{Step}_{z,v}(\cdot)$ for $L$ times,*

$$(x^{(1)}, w^{(1)}) = \mathsf{Step}_{z,v}(x^{(0)}, w^{(0)}),\ (x^{(2)}, w^{(2)}) = \mathsf{Step}_{z,v}(x^{(1)}, w^{(1)}), \dots,\ (x^{(L)}, w^{(L)}) = \mathsf{Step}_{z,v}(x^{(L-1)}, w^{(L-1)}) \ ,$$

*and then define the averages*

$$(\widetilde{x}, \widetilde{w}) = \frac{1}{L}\sum_{\ell=1}^{L} (x^{(\ell)}, w^{(\ell)})$$

*and averaged expansion points $(\widetilde{z}, \widetilde{v})$ as*

$$\widetilde{z} = \frac{1}{L}\sum_{\ell=1}^{L} x^{(\ell)} \quad \text{and} \quad \widetilde{v} = \frac{1}{L}\sum_{\ell=1}^{L} \left( \Sigma^{-1}(w^{(\ell-1)} - w^{(\ell)}) + K(2x^{(\ell)} - x^{(\ell-1)}) \right) \ .$$

*Then there exists a constant $C_{\mathsf{iter}} < \infty$, depending only on $C_{\mathsf{matrix}}, C_{\mathsf{cvx}}, C_{\mathsf{Lip}}, C_{\mathsf{grad}}$ (and in particular, not dependent on L), such that*

$$\left\|\left(\begin{array}{c}\widetilde{z}-x^\star_{z,v}\\ \widetilde{v}-Kx^\star_{z,v}\end{array}\right)\right\|_\Theta \leq C_{\mathsf{iter}}\left(\frac{1}{\sqrt{L}}\left\|\left(\begin{array}{c}x^{(0)}-x^\star_{z,v}\\ w^{(0)}-w^\star_{z,v}\end{array}\right)\right\|_2+\tau R\right).$$

*and*

$$\left\|\left(\begin{array}{c}\widetilde{x}-x^\star_{z,v}\\ \widetilde{w}-w^\star_{z,v}\end{array}\right)\right\|_2 \leq C_{\mathsf{Lip}}+C_{\mathsf{iter}}\left(\frac{1}{\sqrt{L}}\left\|\left(\begin{array}{c}x^{(0)}-x^\star_{z,v}\\ w^{(0)}-w^\star_{z,v}\end{array}\right)\right\|_2+\tau R\right)$$

### 6.2.3 Proof of Theorem 2

We will assume for simplicity that the expansion point is initialized with some $z_0$ satisfying $\|z_0\|_{\mathsf{restrict}} \leq R$; if this is not the case at step $t=0$, then our results can be easily adjusted since we will have $z_1 = x_1 \in \mathsf{dom}(\mathsf{G}_{z_0})$ and therefore $\|z_1\|_{\mathsf{restrict}} \leq R$, so we can simply shift our calculations by one time point.

First, choose any $\delta$ such that

$$0 < \delta < (1-C_{\mathsf{contr}})^{-2}-1$$

and any $L_{\min} > (C_{\mathsf{iter}})^2$.

To prove the desired result, we will prove that

$$\left\|\left(\begin{array}{c}z_{t+1}-x^\star\\ v_{t+1}-Kx^\star\end{array}\right)\right\|_\Theta \leq \frac{C_1}{\sqrt{L'_{t+1}}}+C_2\tau R \tag{35}$$

and

$$\left\|\left(\begin{array}{c}x_{t+1}-x^\star_{z_{t+1},v_{t+1}}\\ w_{t+1}-w^\star_{z_{t+1},v_{t+1}}\end{array}\right)\right\|_2 \leq C_3 \tag{36}$$

for all $t \geq 0$, where

$$C_1 := \max\left\{\sqrt{L'_1}\cdot\left\|\left(\begin{array}{c}z_1-x^\star\\ v_1-Kx^\star\end{array}\right)\right\|_\Theta, \frac{C_{\mathsf{iter}}C_3}{1-(1-C_{\mathsf{contr}})\sqrt{1+\delta}}\right\}$$

and

$$C_2 := \frac{C_{\mathsf{iter}}+C_{\mathsf{excess}}}{C_{\mathsf{contr}}}$$

and

$$C_2 := \max\left\{\left\|\left(\begin{array}{c}x_1-x^\star_{z_1,v_1}\\ w_1-w^\star_{z_1,v_1}\end{array}\right)\right\|_2, \frac{3C_{\mathsf{Lip}}+C_1\cdot\frac{2C_{\mathsf{excess}}}{\sqrt{L_{\min}}}+(C_{\mathsf{iter}}+2C_{\mathsf{excess}}(C_2+1))\tau R}{1-C_{\mathsf{iter}}/\sqrt{L_{\min}}}\right\}.$$

(Note that, if we choose a sufficiently large value for $L_{\min}$, then we can find finite $C_1, C_2, C_3$ satisfying these mutual constraints.)

Assuming that these bounds hold, we then have

$$\|x_{t+1}-x^\star\|_2 = \|z_{t+1}-x^\star\|_2 \leq \left\|\Theta^{-1}\right\|\cdot\left(\frac{C_1}{\sqrt{L'_{t+1}}}+C_2\tau R\right),$$

where the first step holds by definition of $z_{t+1}$; this proves the desired theorem with $C_{\mathsf{converge}} := \left\|\Theta^{-1}\right\|\cdot\max\{C_1, C_2\}$.

Now we prove (35) and (36) by induction. For $t=0$, both statements are true trivially by choice of $C_1, C_2$. Now we will assume that the statements are true for all $t = 0, \ldots, m-1$ and will prove that they hold for $t = m$.

First, for (35), we have

$$\left\|\begin{pmatrix} z_{m+1}-x^\star \\ v_{m+1}-Kx^\star \end{pmatrix}\right\|_\Theta \leq \left\|\begin{pmatrix} z_{m+1}-x^\star_{z_m,v_m} \\ v_{m+1}-Kx^\star_{z_m,v_m} \end{pmatrix}\right\|_\Theta + \left\|\begin{pmatrix} x^\star_m-x^\star \\ Kx^\star_{z_m,v_m}-Kx^\star \end{pmatrix}\right\|_\Theta \quad \text{by the triangle inequality}$$
$$\leq \left\|\begin{pmatrix} z_{m+1}-x^\star_{z_m,v_m} \\ v_{m+1}-Kx^\star_{z_m,v_m} \end{pmatrix}\right\|_\Theta + (1-C_{\mathsf{contr}})\cdot\left\|\begin{pmatrix} z_m-x^\star \\ v_m-Kx^\star \end{pmatrix}\right\|_\Theta + C_{\mathsf{excess}}\tau R \quad \text{by Lemma 1}$$
$$\leq C_{\mathsf{iter}}\left(\frac{1}{\sqrt{L_{m+1}}}\left\|\begin{pmatrix} x_m-x^\star_{z_m,v_m} \\ w_m-w^\star_{z_m,v_m} \end{pmatrix}\right\|_2+\tau R\right) + (1-C_{\mathsf{contr}})\cdot\left\|\begin{pmatrix} z_m-x^\star \\ v_m-Kx^\star \end{pmatrix}\right\|_\Theta + C_{\mathsf{excess}}\tau R \quad \text{by Lemma 2}$$
$$\leq C_{\mathsf{iter}}\left(\frac{1}{\sqrt{L_{m+1}}}C_3+\tau R\right) + (1-C_{\mathsf{contr}})\cdot\left(\frac{C_1}{\sqrt{L'_m}}+C_2\tau R\right) + C_{\mathsf{excess}}\tau R \quad \text{by (35) and (36) applied with } t=m-1$$
$$\leq C_{\mathsf{iter}}\left(\frac{1}{\sqrt{L'_{m+1}}}C_3+\tau R\right) + (1-C_{\mathsf{contr}})\cdot\left(\frac{C_1\sqrt{1+\delta}}{\sqrt{L'_{m+1}}}+C_2\tau R\right) + C_{\mathsf{excess}}\tau R \quad \text{since } L'_{m+1}\leq L_{m+1},(1+\delta)L'_m$$
$$\leq \frac{C_1}{\sqrt{L'_{m+1}}}+C_2\tau R\ ,$$

where the last step holds by definition of $C_1, C_2$. This concludes the proof of (35).

Next, we turn to the proof of (36). For any $t\geq 1$, applying Lemma 1 with $z=z_t$,

$$\left\|\begin{pmatrix} x^\star_{z_t,v_t}-x^\star \\ w^\star_{z_t,v_t}-w^\star \end{pmatrix}\right\|_2 \leq C_{\mathsf{Lip}}+C_{\mathsf{excess}}\left(\left\|\begin{pmatrix} z_t-x^\star \\ v_t-Kx^\star \end{pmatrix}\right\|_\Theta+\tau R\right) \quad \text{by Lemma 1}$$
$$= C_{\mathsf{Lip}}+C_{\mathsf{excess}}\left(\frac{C_1}{\sqrt{L'_t}}+(C_2+1)\tau R\right) \quad \text{by (35) with } t-1 \text{ in place of } t\ ,$$

and also, we have

$$\left\|\begin{pmatrix} x_{m+1}-x^\star_{z_m,v_m} \\ w_{m+1}-w^\star_{z_m,v_m} \end{pmatrix}\right\|_2 \leq C_{\mathsf{Lip}}+C_{\mathsf{iter}}\left(\frac{1}{\sqrt{L_{m+1}}}\left\|\begin{pmatrix} x_m-x^\star_{z_m,v_m} \\ w_m-w^\star_{z_m,v_m} \end{pmatrix}\right\|_2+\tau R\right) \quad \text{by Lemma 2}$$
$$\leq C_{\mathsf{Lip}}+C_{\mathsf{iter}}\left(\frac{1}{\sqrt{L_{m+1}}}C_3+\tau R\right) \quad \text{by (36) applied with } t=m-1.$$

Therefore, combining these calculations, since $L'_m, L'_{m+1}, L_{m+1}\geq L_{\min}$,

$$\left\|\begin{pmatrix} x_{m+1}-x^\star_{z_{m+1},v_{m+1}} \\ w_{m+1}-w^\star_{z_{m+1},v_{m+1}} \end{pmatrix}\right\|_2 \leq \left\|\begin{pmatrix} x_{m+1}-x^\star_{z_m,v_m} \\ w_{m+1}-w^\star_{z_m,v_m} \end{pmatrix}\right\|_2 + \left\|\begin{pmatrix} x^\star_{z_m,v_m}-x^\star \\ w^\star_{z_m,v_m}-w^\star \end{pmatrix}\right\|_2 + \left\|\begin{pmatrix} x^\star_{z_{m+1},v_{m+1}}-x^\star \\ w^\star_{z_{m+1},v_{m+1}}-w^\star \end{pmatrix}\right\|_2$$
$$\leq 3C_{\mathsf{Lip}}+C_{\mathsf{iter}}\left(\frac{1}{\sqrt{L_{\min}}}C_3+\tau R\right)+2C_{\mathsf{excess}}\left(\frac{C_1}{\sqrt{L_{\min}}}+(C_2+1)\tau R\right)$$
$$\leq C_3\ ,$$

by definition of $C_3$. This proves the desired bounds (35) and (36) for $t=m$, and thus we have proved the theorem.

# 7 Discussion

We have developed a primal/dual algorithm for minimizing composite objective functions of the form $\mathsf{F}(Kx)+\mathsf{G}(x)$, which is able to handle nondifferentiability and nonconvexity (even strong concavity) in each individual term, beyond what is possible with many existing approaches based on alternating minimization or proximal gradient methods. The key step of the MOCCA algorithm is the careful choice of local convex approximations to $\mathsf{F}$ and $\mathsf{G}$ at each step, which respects the mirroring between the primal and dual variables of the algorithm. Our method allows for accurate and efficient optimization for a range of problems arising in high-dimensional statistics, such as nonconvex total variation penalties (which reduce the bias caused by shrinkage, when compared to using a convex total variation norm), as well as inverse problems in computed tomography (CT) imaging.

Our present theoretical results give a convergence guarantee, in the case that the overall objective function is approximately convex, for a more stable form of the MOCCA algorithm. In future work, we hope to better understand the relative performance of the various forms of the algorithm, and to find a tighter characterization of the convergence behavior of the algorithm. It would also be interesting to consider a more general form of objective function, $\mathsf{F}(z) + \mathsf{G}(x)$ where $\mathsf{F}, \mathsf{G}$ are nonconvex and nondifferentiable, and where instead of the linear constraint $v = Kx$, the variables $v$ and $x$ are linked via a nonlinear map; such an extension would greatly increase the range of applications of the method.

## Acknowledgements

This work was partially supported by NIH grants CA158446, CA182264, and EB018102. The contents of this article are solely the responsibility of the authors and do not necessarily represent the official views of the National Institutes of Health. The authors are grateful to Wooseok Ha for discussions on the connections with existing algorithms, and to collaborators Taly Gilat Schmidt and Xiaochuan Pan for discussions on the application to CT imaging.

# A Proofs of lemmas

## A.1 Proof of Lemma 1

By definition of $x^\star_{z,v}$,

$$0 \in K^\top \partial \mathsf{F}_v(Kx^\star_{z,v}) + \partial \mathsf{G}_z(x^\star_{z,v}) \tag{37}$$

and since $(x^\star, w^\star)$ is a critical point of the original objective function,

$$0 \in K^\top \partial \mathsf{F}_{Kx^\star}(Kx^\star) + \partial \mathsf{G}_{x^\star}(x^\star)\ . \tag{38}$$

Since $(\mathsf{F}_v - \mathsf{F}_{Kx^\star}) = (\mathsf{F} - \mathsf{F}_{Kx^\star}) - (\mathsf{F} - \mathsf{F}_v)$ and $(\mathsf{G}_z - \mathsf{G}_{x^\star}) = (\mathsf{G} - \mathsf{G}_{x^\star}) - (\mathsf{G} - \mathsf{G}_z)$ are differentiable, we can rewrite (38) as

$$0 \in K^\top \partial \mathsf{F}_v(Kx^\star) + K^\top \nabla(\mathsf{F}-\mathsf{F}_v)(Kx^\star) - K^\top \nabla(\mathsf{F}-\mathsf{F}_{Kx^\star})(Kx^\star) + \partial \mathsf{G}_z(x^\star) + \nabla(\mathsf{G}-\mathsf{G}_z)(x^\star) - \nabla(\mathsf{G}-\mathsf{G}_{x^\star})(x^\star).$$

By the first-order conditions (4), we know that $\nabla(\mathsf{F} - \mathsf{F}_{Kx^\star})(Kx^\star) = 0$ and $\nabla(\mathsf{G} - \mathsf{G}_{x^\star})(x^\star) = 0$, so this reduces to

$$0 \in K^\top \partial \mathsf{F}_v(Kx^\star) + K^\top \nabla(\mathsf{F} - \mathsf{F}_v)(Kx^\star) + \partial \mathsf{G}_z(x^\star) + \nabla(\mathsf{G} - \mathsf{G}_z)(x^\star). \tag{39}$$

We also see that $\|x^\star\|_{\mathsf{restrict}}, \left\|x^\star_{z,v}\right\|_{\mathsf{restrict}} \le R$, since $x^\star, x^\star_{z,v}$ must lie in $\mathsf{dom}(\mathsf{G}) = \mathsf{dom}(\mathsf{G}_z)$.

Then we have

$$\begin{aligned}
&\left\|x^\star_{z,v} - x^\star\right\|^2_{\Theta_\mathsf{G}} + \left\|Kx^\star_{z,v} - Kx^\star\right\|^2_{\Theta_\mathsf{F}} + C_{\mathsf{cvx}} \left\|x^\star_{z,v} - x^\star\right\|^2_2 - 2\tau^2 \cdot \left\|x^\star_{z,v} - x^\star\right\|^2_{\mathsf{restrict}}\\
&\le \left\|Kx^\star_{z,v} - Kx^\star\right\|^2_{\Lambda_\mathsf{F}} + \left\|x^\star_{z,v} - x^\star\right\|^2_{\Lambda_\mathsf{G}} - \tau^2 \cdot \left\|x^\star_{z,v} - x^\star\right\|^2_{\mathsf{restrict}} \quad \text{by Assumption 2}\\
&\le \langle Kx^\star_{z,v} - Kx^\star, \partial \mathsf{F}_u(Kx^\star_{z,v}) - \partial \mathsf{F}_u(Kx^\star)\rangle + \langle x^\star_{z,v} - x^\star, \partial \mathsf{G}_z(x^\star_{z,v}) - \partial \mathsf{G}_z(x^\star)\rangle \quad \text{by Assumption 2}\\
&= \langle Kx^\star_{z,v} - Kx^\star, \nabla(\mathsf{F} - \mathsf{F}_v)(Kx^\star)\rangle + \langle x^\star_{z,v} - x^\star, \nabla(\mathsf{G} - \mathsf{G}_z)(x^\star)\rangle \quad \text{by (37) and (39)}\\
&\le \frac{1}{2}\left\|Kx^\star_{z,v} - Kx^\star\right\|^2_{\Theta_\mathsf{F}} + \frac{1}{2}\|v - Kx^\star\|^2_{\Theta_\mathsf{F}} + \frac{1}{2}\left\|x^\star_{z,v} - x^\star\right\|^2_{\Theta_\mathsf{G}} + \frac{1}{2}\|z - x^\star\|^2_{\Theta_\mathsf{G}}\\
&\qquad + \tau^2\left(\left\|x^\star_{z,v} - x^\star\right\|^2_{\mathsf{restrict}} + \|z - x^\star\|^2_{\mathsf{restrict}}\right) \quad \text{by Assumption 2}\ .
\end{aligned}$$

Next, recall that $\left\|x^\star_{z,v}\right\|_{\textsf{restrict}}, \|x^\star\|_{\textsf{restrict}} \leq R$ from before and $\|z\|_{\textsf{restrict}} \leq R$ by assumption. After rearranging terms, this gives

$$\frac{1}{2}\left\|x^\star_{z,v} - x^\star\right\|^2_{\Theta_{\mathsf{G}}} + \frac{1}{2}\left\|Kx^\star_{z,v} - Kx^\star\right\|^2_{\Theta_{\mathsf{F}}} + C_{\textsf{cvx}}\left\|x^\star_{z,v} - x^\star\right\|^2_2 \leq \frac{1}{2}\|v - Kx^\star\|^2_{\Theta_{\mathsf{F}}} + \frac{1}{2}\|z - x^\star\|^2_{\Theta_{\mathsf{G}}} + 16\tau^2 R^2 \ . \tag{40}$$

Now, using the definition of $\Theta$, we see that

$$\begin{pmatrix} \Theta_{\mathsf{G}} & 0 \\ 0 & \Theta_{\mathsf{F}} \end{pmatrix} \preceq (1 - C_{\textsf{contr}})^2 \cdot \Theta$$

where

$$C_{\textsf{contr}} := \min\left\{\frac{1}{2}, \frac{C_{\textsf{cvx}}}{8(C_{\textsf{matrix}})^3}\right\} > 0 \ .$$

We then get

$$\begin{aligned}
&(1 - C_{\textsf{contr}})^2 \cdot \frac{1}{2}\left\|\begin{pmatrix} z - x^\star \\ v - Kx^\star \end{pmatrix}\right\|^2_{\Theta} + 16\tau^2 R^2 \\
&\geq \frac{1}{2}\|z - x^\star\|^2_{\Theta_{\mathsf{G}}} + \frac{1}{2}\|v - Kx^\star\|^2_{\Theta_{\mathsf{F}}} + 16\tau^2 R^2 \\
&\geq \frac{1}{2}\left\|x^\star_{z,v} - x^\star\right\|^2_{\Theta_{\mathsf{G}}} + \frac{1}{2}\left\|Kx^\star_{z,v} - Kx^\star\right\|^2_{\Theta_{\mathsf{F}}} + C_{\textsf{cvx}}\left\|x^\star_{z,v} - x^\star\right\|^2_2 \quad \text{from (40)} \\
&\geq \frac{1}{2}\left\|x^\star_{z,v} - x^\star\right\|^2_{\Theta_{\mathsf{G}}} + \frac{1}{2}\left\|Kx^\star_{z,v} - Kx^\star\right\|^2_{\Theta_{\mathsf{F}}} + \frac{C_{\textsf{cvx}}}{2}\left\|x^\star_{z,v} - x^\star\right\|^2_2 + \frac{C_{\textsf{cvx}}}{2(C_{\textsf{matrix}})^2}\left\|Kx^\star_{z,v} - Kx^\star\right\|^2_2 \\
&\geq \frac{1}{2}\left\|x^\star_{z,v} - x^\star\right\|^2_{\Theta_{\mathsf{G}} + C_{\textsf{cvx}}/(C_{\textsf{matrix}})^2 \cdot \mathbf{I}} + \frac{1}{2}\left\|Kx^\star_{z,v} - Kx^\star\right\|^2_{\Theta_{\mathsf{F}} + C_{\textsf{cvx}}/(C_{\textsf{matrix}})^2 \cdot \mathbf{I}} \\
&= \frac{1}{2}\left\|\begin{pmatrix} x^\star_{z,v} - x^\star \\ Kx^\star_{z,v} - Kx^\star \end{pmatrix}\right\|^2_{\Theta} \ .
\end{aligned}$$

Rearranging terms, we obtain

$$\left\|\begin{pmatrix} x^\star_{z,v} - x^\star \\ Kx^\star_{z,v} - Kx^\star \end{pmatrix}\right\|_{\Theta} \leq (1 - C_{\textsf{contr}})\left\|\begin{pmatrix} z - x^\star \\ v - Kx^\star \end{pmatrix}\right\|_{\Theta} + \sqrt{32}\tau R \ .$$

Next, by definition, we also have

$$w^\star_{z,v} \in \partial\mathsf{F}_v(Kx^\star_{z,v}) \text{ and } w^\star \in \partial\mathsf{F}_{Kx^\star}(Kx^\star) \ .$$

This second expression can be rewritten as

$$w^\star \in \partial\mathsf{F}_v(Kx^\star) + \nabla(\mathsf{F} - \mathsf{F}_v)(Kx^\star) - \nabla(\mathsf{F} - \mathsf{F}_{Kx^\star})(Kx^\star) = \partial\mathsf{F}_v(Kx^\star) + \nabla(\mathsf{F} - \mathsf{F}_v)(Kx^\star) \ ,$$

where the last step uses the first-order condition (4) to see that $\nabla(\mathsf{F} - \mathsf{F}_{Kx^\star})(Kx^\star) = 0$. Therefore,

$$\begin{aligned}
\left\|w^\star_{z,v} - w^\star\right\|_2 &\leq \left\|\partial\mathsf{F}_v(Kx^\star) - \partial\mathsf{F}_v(Kx^\star_{z,v})\right\|_2 + \|\nabla(\mathsf{F} - \mathsf{F}_v)(Kx^\star)\|_2 \\
&\leq C_{\textsf{Lip}} + C_{\textsf{grad}}\left\|Kx^\star_{z,v} - Kx^\star\right\|_2 + \sqrt{\|\Theta_{\mathsf{F}}\|}\,\|v - Kx^\star\|_{\Theta_{\mathsf{F}}} \quad \text{by Assumption 2} \\
&\leq C_{\textsf{Lip}} + C_{\textsf{grad}}C_{\textsf{matrix}}\left\|x^\star_{z,v} - x^\star\right\|_2 + \sqrt{C_{\textsf{matrix}}}\,\|v - Kx^\star\|_{\Theta_{\mathsf{F}}} \ ,
\end{aligned}$$

and so, by the work above,

$$\begin{aligned}\left\|\left(\begin{array}{c} x^\star_{z,v}-x^\star \\ w^\star_{z,v}-w^\star \end{array}\right)\right\|_2 &\leq C_{\mathsf{Lip}} + (1+C_{\mathsf{grad}}C_{\mathsf{matrix}})\left\|x^\star_{z,v}-x^\star\right\|_2 + \sqrt{C_{\mathsf{matrix}}}\left\|v-Kx^\star\right\|_{\Theta_{\mathsf{F}}} \\ &\leq C_{\mathsf{Lip}} + (1+C_{\mathsf{grad}}C_{\mathsf{matrix}})\frac{\sqrt{(1-C_{\mathsf{contr}})^2\cdot\frac{1}{2}\left\|\left(\begin{array}{c} z-x^\star \\ v-Kx^\star \end{array}\right)\right\|^2_\Theta + 16\tau^2R^2}}{\sqrt{C_{\mathsf{cvx}}}} + \sqrt{C_{\mathsf{matrix}}}\left\|v-Kx^\star\right\|_{\Theta_{\mathsf{F}}} \\ &\leq C_{\mathsf{Lip}} + (1+C_{\mathsf{grad}}C_{\mathsf{matrix}})\frac{\frac{1}{\sqrt{2}}\left\|\left(\begin{array}{c} z-x^\star \\ v-Kx^\star \end{array}\right)\right\|_\Theta + 4\tau R}{\sqrt{C_{\mathsf{cvx}}}} + \sqrt{C_{\mathsf{matrix}}}\left\|v-Kx^\star\right\|_{\Theta_{\mathsf{F}}} \\ &\leq C_{\mathsf{Lip}} + (1+C_{\mathsf{grad}}C_{\mathsf{matrix}})\frac{\frac{1}{\sqrt{2}}\left\|\left(\begin{array}{c} z-x^\star \\ v-Kx^\star \end{array}\right)\right\|_\Theta + 4\tau R}{\sqrt{C_{\mathsf{cvx}}}} + \sqrt{C_{\mathsf{matrix}}}\left\|\left(\begin{array}{c} z-x^\star \\ v-Kx^\star \end{array}\right)\right\|_\Theta ,\end{aligned}$$

by definition of $\Theta$. Setting

$$C_{\mathsf{excess}} = \max\left\{32, \sqrt{C_{\mathsf{matrix}}} + 4\cdot\frac{1+C_{\mathsf{grad}}C_{\mathsf{matrix}}}{\sqrt{C_{\mathsf{cvx}}}},\right\} ,$$

this is sufficient to prove the lemma.

## A.2 Proof of Lemma 2

First we state and prove a supporting lemma which considers only a single step of the "inner loop".

**Lemma 3.** *There exists a constant $C_{\mathsf{monotone}} > 0$, which depends only on $C_{\mathsf{matrix}}, C_{\mathsf{cvx}}, C_{\mathsf{Lip}}, C_{\mathsf{grad}}$, such that, for any $x, y, z$ with $\|x\|_{\mathsf{restrict}} \leq R$, if $(x', w') = \mathsf{Step}_{z,v}(x, w)$, then*

$$C_{\mathsf{monotone}}\left(\left\|x'-x^\star_{z,v}\right\|^2_2 + \left\|\left(\begin{array}{c} x-x' \\ w-w' \end{array}\right)\right\|^2_2\right) \leq \left\|\left(\begin{array}{c} x-x^\star_{z,v} \\ w-w^\star_{z,v} \end{array}\right)\right\|^2_M - \left\|\left(\begin{array}{c} x'-x^\star_{z,v} \\ w'-w^\star_{z,v} \end{array}\right)\right\|^2_M + \tau^2R^2 .$$

*Proof of Lemma 3.* By definition of the $(x, w)$ update step, we have

$$\mathrm{T}^{-1}(x-x') - K^\top w \in \partial\mathsf{G}_z(x') \tag{41}$$

and

$$\Sigma^{-1}(w-w') + K(2x'-x) \in \partial\mathsf{F}^*_v(w') ,$$

and from (32) we know that this last expression implies

$$w' \in \partial\mathsf{F}_v\left(\Sigma^{-1}(w-w') + K(2x'-x)\right) . \tag{42}$$

Similarly, by definition of $(x^\star_{z,v}, w^\star_{z,v})$, we also have

$$-K^\top w^\star_{z,v} \in \partial\mathsf{G}_z(x^\star_{z,v}) \tag{43}$$

and

$$w^\star_{z,v} \in \partial\mathsf{F}_v(Kx^\star_{z,v}) . \tag{44}$$

Therefore, combining these expressions, we have

$$\left\langle \begin{pmatrix} x' - x^\star_{z,v} \\ \Sigma^{-1}(w - w') + K(2x' - x) - Kx^\star_{z,v} \end{pmatrix}, \begin{pmatrix} \partial\mathsf{G}_z(x') - \partial\mathsf{G}_z(x^\star_{z,v}) \\ \partial\mathsf{F}_v\left(K(2x' - x) + \Sigma^{-1}(w - w')\right) - \partial\mathsf{F}_v(Kx^\star_{z,v}) \end{pmatrix} \right\rangle$$
$$\ni \left\langle \begin{pmatrix} x' - x^\star_{z,v} \\ \Sigma^{-1}(w - w') + K(2x' - x) - Kx^\star_{z,v} \end{pmatrix}, \begin{pmatrix} \mathrm{T}^{-1}(x - x') - K^\top(w - w^\star_{z,v}) \\ w' - w^\star_{z,v} \end{pmatrix} \right\rangle$$
$$= \left\langle \begin{pmatrix} x' - x^\star_{z,v} \\ w' - w^\star_{z,v} \end{pmatrix}, \begin{pmatrix} \mathrm{T}^{-1}(x - x') - K^\top(w - w') \\ \Sigma^{-1}(w - w') - K(x - x') \end{pmatrix} \right\rangle \quad \text{by reorganizing terms}$$
$$= \begin{pmatrix} x' - x^\star_{z,v} \\ w' - w^\star_{z,v} \end{pmatrix}^\top M \begin{pmatrix} x - x' \\ w - w' \end{pmatrix}. \tag{45}$$

As in [7] we can calculate

$$\left| \begin{pmatrix} x' - x^\star_{z,v} \\ w' - w^\star_{z,v} \end{pmatrix}^\top M \begin{pmatrix} x - x' \\ w - w' \end{pmatrix} \right| = \frac{1}{2} \left\| \begin{pmatrix} x - x^\star_{z,v} \\ w - w^\star_{z,v} \end{pmatrix} \right\|^2_M - \frac{1}{2} \left\| \begin{pmatrix} x' - x^\star_{z,v} \\ w' - w^\star_{z,v} \end{pmatrix} \right\|^2_M - \frac{1}{2} \left\| \begin{pmatrix} x - x' \\ w - w' \end{pmatrix} \right\|^2_M .$$

On the other hand, by the convexity of $\mathsf{F}_v$ and $\mathsf{G}_z$ as stated in Assumption 2, we have

$$\left\langle \begin{pmatrix} x' - x^\star_{z,v} \\ \Sigma^{-1}(w - w') + K(2x' - x) - Kx^\star_{z,v} \end{pmatrix}, \begin{pmatrix} \partial\mathsf{G}_z(x') - \partial\mathsf{G}_z(x^\star_{z,v}) \\ \partial\mathsf{F}_v\left(K(2x' - x) + \Sigma^{-1}(w - w')\right) - \partial\mathsf{F}_v(Kx^\star_{z,v}) \end{pmatrix} \right\rangle$$
$$\geq \left\| x' - x^\star_{z,v} \right\|^2_{\Lambda_{\mathsf{G}}} + \left\| \Sigma^{-1}(w - w') + K(2x' - x) - Kx^\star_{z,v} \right\|^2_{\Lambda_{\mathsf{F}}} - \tau^2 \left\| x' - x^\star_{z,v} \right\|^2_{\mathsf{restrict}}$$
$$\geq \left\| x' - x^\star_{z,v} \right\|^2_{\Lambda_{\mathsf{G}}} + \frac{1}{2} \left\| Kx' - Kx^\star_{z,v} \right\|^2_{\Lambda_{\mathsf{F}}} - \left\| \Sigma^{-1}(w - w') - K(x - x') \right\|^2_{\Lambda_{\mathsf{F}}} - \tau^2 \left\| x' - x^\star_{z,v} \right\|^2_{\mathsf{restrict}}$$
$$\geq \frac{C_{\mathsf{cvx}}}{2} \left\| x' - x^\star_{z,v} \right\|^2_2 - \frac{3\tau^2}{2} \left\| x' - x^\star_{z,v} \right\|^2_{\mathsf{restrict}} - \left\| \Sigma^{-1}(w - w') - K(x - x') \right\|^2_{\Lambda_{\mathsf{F}}}$$
$$\geq \frac{C_{\mathsf{cvx}}}{2} \left\| x' - x^\star_{z,v} \right\|^2_2 - \frac{3\tau^2}{2} \left\| x' - x^\star_{z,v} \right\|^2_{\mathsf{restrict}} - \left\| \begin{pmatrix} x - x' \\ w - w' \end{pmatrix} \right\|^2_M \cdot \left\| M^{-1} \right\| \left\| \Lambda_{\mathsf{F}} \right\| \left( \left\| \Sigma^{-1} \right\| + \| K \| \right)$$
$$\geq \frac{C_{\mathsf{cvx}}}{2} \left\| x' - x^\star_{z,v} \right\|^2_2 - 6\tau^2 R^2 - \left\| \begin{pmatrix} x - x' \\ w - w' \end{pmatrix} \right\|^2_M \cdot 2(C_{\mathsf{matrix}})^3 ,$$

where the last step holds because $\|x'\|_{\mathsf{restrict}}, \left\|x^\star_{z,v}\right\|_{\mathsf{restrict}} \leq R$, since $x'$ and $x^\star_{z,v}$ must both lie in $\mathsf{dom}(\mathsf{G}) = \mathsf{dom}(\mathsf{G}_z)$ by their definitions. Now, examining these calculations, we see that the left-hand side must be nonnegative, so we can also write

$$\left\langle \begin{pmatrix} x' - x^\star_{z,v} \\ \Sigma^{-1}(w - w') + K(2x' - x) - Kx^\star_{z,v} \end{pmatrix}, \begin{pmatrix} \partial\mathsf{G}_z(x') - \partial\mathsf{G}_z(x^\star_{z,v}) \\ \partial\mathsf{F}_v\left(K(2x' - x) + \Sigma^{-1}(w - w')\right) - \partial\mathsf{F}_v(Kx^\star_{z,v}) \end{pmatrix} \right\rangle$$
$$\geq c \left( \frac{C_{\mathsf{cvx}}}{2} \left\| x' - x^\star_{z,v} \right\|^2_2 - 6\tau^2 R^2 - \left\| \begin{pmatrix} x - x' \\ w - w' \end{pmatrix} \right\|^2_M \cdot 2(C_{\mathsf{matrix}})^3 \right)$$

for any $c \in [0, 1]$. Choosing $c = \frac{1}{\max\{12, 8C_{\mathsf{matrix}}\}}$, we obtain

$$\left\langle \begin{pmatrix} x' - x^\star_{z,v} \\ \Sigma^{-1}(w - w') + K(2x' - x) - Kx^\star_{z,v} \end{pmatrix}, \begin{pmatrix} \partial\mathsf{G}_z(x') - \partial\mathsf{G}_z(x^\star_{z,v}) \\ \partial\mathsf{F}_v\left(K(2x' - x) + \Sigma^{-1}(w - w')\right) - \partial\mathsf{F}_v(Kx^\star_{z,v}) \end{pmatrix} \right\rangle$$
$$\geq \frac{C_{\mathsf{cvx}}}{2\max\{12, 8C_{\mathsf{matrix}}\}} \left\| x' - x^\star_{z,v} \right\|^2_2 - \frac{1}{2}\tau^2 R^2 - \frac{1}{4} \left\| \begin{pmatrix} x - x' \\ w - w' \end{pmatrix} \right\|^2_M$$

Combining all our work, then,

$$\frac{C_{\mathsf{cvx}}}{4\max\{12, 8C_{\mathsf{matrix}}\}} \left\| x' - x^\star_{z,v} \right\|^2_2 \leq \left\| \begin{pmatrix} x - x^\star_{z,v} \\ w - w^\star_{z,v} \end{pmatrix} \right\|^2_M - \left\| \begin{pmatrix} x' - x^\star_{z,v} \\ w' - w^\star_{z,v} \end{pmatrix} \right\|^2_M - \frac{1}{2} \left\| \begin{pmatrix} x - x' \\ w - w' \end{pmatrix} \right\|^2_M + \tau^2 R^2 .$$

Setting $C_{\textsf{monotone}} = \min\left\{\frac{1}{2C_{\textsf{matrix}}}, \frac{C_{\textsf{cvx}}}{4\max\{12, 8C_{\textsf{matrix}}\}}\right\}$, we have proved the lemma. □

Now we turn to the proof of Lemma 2. By Lemma 3, for each $\ell = 1, \dots, L$, we have

$$C_{\textsf{monotone}}\left(\left\|x^{(\ell)} - x^\star_{z,v}\right\|_2^2 + \left\|\begin{pmatrix} x^{(\ell-1)} - x^{(\ell)} \\ w^{(\ell-1)} - w^{(\ell)} \end{pmatrix}\right\|_2^2\right) \leq \left\|\begin{pmatrix} x^{(\ell-1)} - x^\star_{z,v} \\ w^{(\ell-1)} - w^\star_{z,v} \end{pmatrix}\right\|_M^2 - \left\|\begin{pmatrix} x^{(\ell)} - x^\star_{z,v} \\ w^{(\ell)} - w^\star_{z,v} \end{pmatrix}\right\|_M^2 + \tau^2 R^2 \ .$$

Summing this inequality over $\ell = 1, \dots, L$, and taking a telescoping sum on the right-hand side, we have

$$C_{\textsf{monotone}}\sum_{\ell=1}^L\left(\left\|x^{(\ell)} - x^\star_{z,v}\right\|_2^2 + \left\|\begin{pmatrix} x^{(\ell-1)} - x^{(\ell)} \\ w^{(\ell-1)} - w^{(\ell)} \end{pmatrix}\right\|_2^2\right) \leq \left\|\begin{pmatrix} x^{(0)} - x^\star_{z,v} \\ w^{(0)} - w^\star_{z,v} \end{pmatrix}\right\|_M^2 + L\tau^2 R^2 \ .$$

Next, by convexity of $w \mapsto \|w\|_2^2$, we have

$$\left\|\widetilde{x} - x^\star_{z,v}\right\|_2^2 \leq \frac{1}{L}\sum_{\ell=1}^L\left\|x^{(\ell)} - x^\star_{z,v}\right\|_2^2$$

and

$$\left\|\frac{1}{L}\begin{pmatrix} x^{(0)} - x^{(L)} \\ w^{(0)} - w^{(L)} \end{pmatrix}\right\|_2^2 \leq \frac{1}{L}\sum_{\ell=1}^L\left\|\begin{pmatrix} x^{(\ell-1)} - x^{(\ell)} \\ w^{(\ell-1)} - w^{(\ell)} \end{pmatrix}\right\|_2^2 \ .$$

So,

$$C_{\textsf{monotone}}\left(\left\|\widetilde{x} - x^\star_{z,v}\right\|_2^2 + \left\|\frac{1}{L}\begin{pmatrix} x^{(0)} - x^{(L)} \\ w^{(0)} - w^{(L)} \end{pmatrix}\right\|_2^2\right) \leq \frac{1}{L}\left\|\begin{pmatrix} x^{(0)} - x^\star_{z,v} \\ w^{(0)} - w^\star_{z,v} \end{pmatrix}\right\|_M^2 + \tau^2 R^2 \ .$$

Next, by definition of $\widetilde{z}$, we can write

$$\left\|\begin{pmatrix} \widetilde{z}_{\textsf{G}} - x^\star_{z,v} \\ \widetilde{z}_{\textsf{F}} - Kx^\star_{z,v} \end{pmatrix}\right\|_\Theta \leq \left\|\begin{pmatrix} \widetilde{x} - x^\star_{z,v} \\ K\widetilde{x} - Kx^\star_{z,v} \end{pmatrix}\right\|_\Theta + \left\|\frac{1}{L}\begin{pmatrix} 0 \\ \Sigma^{-1}(w^{(0)} - w^{(L)}) + K(x^{(L)} - x^{(0)}) \end{pmatrix}\right\|_\Theta$$

and so clearly,

$$\left\|\begin{pmatrix} \widetilde{z}_{\textsf{G}} - x^\star_{z,v} \\ \widetilde{z}_{\textsf{F}} - Kx^\star_{z,v} \end{pmatrix}\right\|_\Theta^2 \leq \frac{\|\Theta\|\,(1 + 2\,\|K\| + \|\Sigma^{-1}\|)^2}{C_{\textsf{monotone}}}\left(\frac{1}{L}\left\|\begin{pmatrix} x^{(0)} - x^\star_{z,v} \\ w^{(0)} - w^\star_{z,v} \end{pmatrix}\right\|_M^2 + \tau^2 R^2\right) \ .$$

Finally, for each $\ell$, by definition of the step,

$$w^{(\ell)} \in \partial\textsf{F}_v\left(\Sigma^{-1}(w^{(\ell-1)} - w^{(\ell)}) + K\bar{x}^{(\ell)}\right)$$

while

$$w^\star_{z,v} \in \partial\textsf{F}_v(Kx^\star_{z,v}) \ .$$

Therefore, by Assumption 2,

$$\left\|w^{(\ell)} - w^\star_{z,v}\right\|_2 \leq C_{\textsf{Lip}} + C_{\textsf{grad}}\left\|\Sigma^{-1}(w^{(\ell-1)} - w^{(\ell)}) + K\bar{x}^{(\ell)} - Kx^\star_{z,v}\right\|_2 \ .$$

By convexity, then,

$$\begin{aligned}\left\|\widetilde{w}-w^\star_{z,v}\right\|_2 &\leq \frac{1}{L}\sum_{\ell=1}^{L}\left(C_{\mathsf{Lip}}+C_{\mathsf{grad}}\left\|\Sigma^{-1}(w^{(\ell-1)}-w^{(\ell)})+K\bar{x}^{(\ell)}-Kx^\star_{z,v}\right\|_2\right)\\ &\leq C_{\mathsf{Lip}}+C_{\mathsf{grad}}(\|K\|+\left\|\Sigma^{-1}\right\|)\left(\frac{1}{L}\sum_{\ell=1}^{L}\left\|x^{(\ell)}-x^\star_{z,v}\right\|_2+\left\|\begin{pmatrix}x^{(\ell-1)}-x^{(\ell)}\\ w^{(\ell-1)}-w^{(\ell)}\end{pmatrix}\right\|_2\right)\\ &\leq C_{\mathsf{Lip}}+C_{\mathsf{grad}}(\|K\|+\left\|\Sigma^{-1}\right\|)\sqrt{\frac{1}{L}\sum_{\ell=1}^{L}\|x^{(\ell)}-x^\star_{z,v}\|_2^2+\left\|\begin{pmatrix}x^{(\ell-1)}-x^{(\ell)}\\ w^{(\ell-1)}-w^{(\ell)}\end{pmatrix}\right\|_2^2}\\ &\leq C_{\mathsf{Lip}}+C_{\mathsf{grad}}(\|K\|+\left\|\Sigma^{-1}\right\|)\sqrt{\frac{1}{C_{\mathsf{monotone}}}\left(\frac{1}{L}\left\|\begin{pmatrix}x^{(0)}-x^\star_{z,v}\\ w^{(0)}-w^\star_{z,v}\end{pmatrix}\right\|_M^2+\tau^2R^2\right)}\ .\end{aligned}$$

Combining everything, this proves that

$$\left\|\begin{pmatrix}\widetilde{x}-x^\star_{z,v}\\ \widetilde{w}-w^\star_{z,v}\end{pmatrix}\right\|_2 \leq C_{\mathsf{Lip}}+C_{\mathsf{grad}}(1+\|K\|+\left\|\Sigma^{-1}\right\|)\left(1+\frac{1+\sqrt{\|M\|}}{\sqrt{C_{\mathsf{monotone}}}}\right)\cdot\left(\frac{1}{\sqrt{L}}\left\|\begin{pmatrix}x^{(0)}-x^\star_{z,v}\\ w^{(0)}-w^\star_{z,v}\end{pmatrix}\right\|_2+\tau R\right)\ .$$

Finally, defining

$$C_{\mathsf{iter}}=\max\left\{\sqrt{\frac{C_{\mathsf{matrix}}(1+3C_{\mathsf{matrix}})^2}{C_{\mathsf{monotone}}}},C_{\mathsf{grad}}(1+2C_{\mathsf{matrix}})\left(1+\frac{1+\sqrt{C_{\mathsf{matrix}}}}{\sqrt{C_{\mathsf{monotone}}}}\right)\right\}$$

and using the definition of $C_{\mathsf{matrix}}$, we have proved the lemma.

# B Connection between MOCCA and ADMM

We now give the details for the connection between MOCCA and the Alternating Directions Method of Multipliers (ADMM) algorithm [4], which we introduced in Section 3.3. We reformulate the basic version of MOCCA, Algorithm 1, as an ADMM type algorithm with preconditioning and with convex approximations to $\mathsf{F}$ and $\mathsf{G}$ (the stable "inner loop" version, Algorithm 2, can also be interpreted as an extension of ADMM, but we do not give details here). The resulting algorithm is given in Algorithm 3. The equivalence between Algorithms 1 and 3 is simply an extension of the connection between the Chambolle-Pock algorithm and a preconditioned ADMM, as shown in [7].

**Algorithm 3** MOCCA algorithm: ADMM version

**Input:** Convex functions $\mathsf{F}_{\mathsf{cvx}},\mathsf{G}_{\mathsf{cvx}}$, differentiable functions $\mathsf{F}_{\mathsf{diff}},\mathsf{G}_{\mathsf{diff}}$, linear operator $K$, positive diagonal step size matrices $\Sigma,\mathrm{T}$, extrapolation parameter $\theta\in[0,1]$.
**Initialize:** Primal variables $x_0\in\mathbb{R}^d$, $u_0\in\mathbb{R}^m$, dual variable $\Delta_0\in\mathbb{R}^m$.
**for** $t=0,1,2,\dots$ **do**
Update all variables:

$$\begin{cases}x_{t+1}=\arg\min_x\left\{\mathsf{G}_{x_t}(x)+\langle x,K^\top\Delta_t\rangle+\frac{1}{2}\left\|Kx-u_t\right\|_\Sigma^2+\frac{1}{2}\left\|x-x_t\right\|^2_{\mathrm{T}^{-1}-K^\top\Sigma K}\right\},\\ \Delta_{t+1}=\Delta_t+\Sigma(Kx_{t+1}-u_t),\\ u_{t+1}=\arg\min_u\left\{\mathsf{F}_{u_t}(u)-\langle\Delta_{t+1},u\rangle+\frac{1}{2}\left\|Kx_{t+1}-u\right\|_\Sigma^2\right\},\end{cases}$$

**until** some convergence criterion is reached.

# References

[1] Rina Foygel Barber, Emil Y. Sidky, Taly Gilat Schmidt, and Xiaochuan Pan. An algorithm for constrained one-step inversion of spectral CT data. *Physics in Medicine and Biology*, 61(10):3784–3818, 2016.

[2] Amir Beck and Marc Teboulle. A fast iterative shrinkage-thresholding algorithm for linear inverse problems. *SIAM Journal on Imaging Sciences*, 2(1):183–202, 2009.

[3] Jérôme Bolte, Shoham Sabach, and Marc Teboulle. Proximal alternating linearized minimization for nonconvex and nonsmooth problems. *Mathematical Programming*, 146(1-2):459–494, 2014.

[4] Stephen Boyd, Neal Parikh, Eric Chu, Borja Peleato, and Jonathan Eckstein. Distributed optimization and statistical learning via the alternating direction method of multipliers. *Foundations and Trends in Machine Learning*, 3(1):1–122, 2011.

[5] Emmanuel J Candès, Michael B Wakin, and Stephen P Boyd. Enhancing sparsity by reweighted $\ell_1$ minimization. *Journal of Fourier Analysis and Applications*, 14(5-6):877–905, 2008.

[6] Antonin Chambolle and Jérôme Darbon. On total variation minimization and surface evolution using parametric maximum flows. *International journal of computer vision*, 84(3):288–307, 2009.

[7] Antonin Chambolle and Thomas Pock. A first-order primal-dual algorithm for convex problems with applications to imaging. *Journal of Mathematical Imaging and Vision*, 40(1):120–145, 2011.

[8] Antonin Chambolle and Thomas Pock. A remark on accelerated block coordinate descent for computing the proximity operators of a sum of convex functions. *SMAI Journal of Computational Mathematics*, 1: 29–54, 2015.

[9] Rick Chartrand. Exact reconstruction of sparse signals via nonconvex minimization. *IEEE Signal Processing Letters*, 14(10):707–710, 2007.

[10] Aditya Chopra and Heng Lian. Total variation, adaptive total variation and nonconvex smoothly clipped absolute deviation penalty for denoising blocky images. *Pattern Recognition*, 43(8):2609–2619, 2010.

[11] Ernie Esser, Xiaoqun Zhang, and Tony Chan. A general framework for a class of first order primal-dual algorithms for TV minimization. *UCLA CAM Report*, pages 09–67, 2009.

[12] Jianqing Fan and Runze Li. Variable selection via nonconcave penalized likelihood and its oracle properties. *Journal of the American Statistical Association*, 96(456):1348–1360, 2001.

[13] Saeed Ghadimi and Guanghui Lan. Accelerated gradient methods for nonconvex nonlinear and stochastic programming. *Mathematical Programming*, 156(1):59–99, 2016.

[14] Bingsheng He and Xiaoming Yuan. Convergence analysis of primal-dual algorithms for a saddle-point problem: from contraction perspective. *SIAM Journal on Imaging Sciences*, 5(1):119–149, 2012.

[15] Mingyi Hong, Zhi-Quan Luo, and Meisam Razaviyayn. Convergence analysis of alternating direction method of multipliers for a family of nonconvex problems. *arXiv preprint arXiv:1410.1390*, 2014.

[16] David R Hunter and Kenneth Lange. Quantile regression via an mm algorithm. *Journal of Computational and Graphical Statistics*, 9(1):60–77, 2000.

[17] Nicholas A Johnson. A dynamic programming algorithm for the fused lasso and l 0-segmentation. *Journal of Computational and Graphical Statistics*, 22(2):246–260, 2013.

[18] Keith Knight and Wenjiang Fu. Asymptotics for lasso-type estimators. *Annals of Statistics*, 28(5):1356–1378, 2000.

[19] Guoyin Li and Ting Kei Pong. Global convergence of splitting methods for nonconvex composite optimization. *arXiv preprint arXiv:1407.0753*, 2014.

[20] Huan Li and Zhouchen Lin. Accelerated proximal gradient methods for nonconvex programming. In *Advances in Neural Information Processing Systems*, pages 379–387, 2015.

[21] Po-Ling Loh and Martin J Wainwright. High-dimensional regression with noisy and missing data: Provable guarantees with non-convexity. In *Advances in Neural Information Processing Systems*, pages 2726–2734, 2011.

[22] Po-Ling Loh and Martin J Wainwright. Regularized M-estimators with nonconvexity: Statistical and algorithmic theory for local optima. In *Advances in Neural Information Processing Systems*, pages 476–484, 2013.

[23] Chengwu Lu and Hua Huang. TV + $TV_2$ regularization with nonconvex sparseness-inducing penalty for image restoration. *Mathematical Problems in Engineering*, 2014:790547, 2014.

[24] Sindri Magnússon, Pradeep Chathuranga Weeraddana, Michael G Rabbat, and Carlo Fischione. On the convergence of alternating direction lagrangian methods for nonconvex structured optimization problems. *arXiv preprint arXiv:1409.8033*, 2014.

[25] MATLAB. *Version 8.6.0 (R2015b)*. The MathWorks Inc., Natick, Massachusetts, 2015.

[26] Sahand Negahban, Bin Yu, Martin J Wainwright, and Pradeep K Ravikumar. A unified framework for high-dimensional analysis of M-estimators with decomposable regularizers. In *Advances in Neural Information Processing Systems*, pages 1348–1356, 2009.

[27] Yurii Nesterov. Gradient methods for minimizing composite objective function, 2007. URL `http://hdl.handle.net/2078.1/5122`.

[28] Peter Ochs, Yunjin Chen, Thomas Brox, and Thomas Pock. iPiano: inertial proximal algorithm for nonconvex optimization. *SIAM Journal on Imaging Sciences*, 7(2):1388–1419, 2014.

[29] Peter Ochs, Alexey Dosovitskiy, Thomas Brox, and Thomas Pock. On iteratively reweighted algorithms for nonsmooth nonconvex optimization in computer vision. *SIAM Journal on Imaging Sciences*, 8:331–372, 2015.

[30] James M Ortega and Werner C Rheinboldt. *Iterative solution of nonlinear equations in several variables*, volume 30. SIAM, 1970.

[31] Ankit Parekh and Ivan W Selesnick. Convex fused lasso denoising with non-convex regularization and its use for pulse detection. *arXiv preprint arXiv:1509.02811*, 2015.

[32] Thomas Pock and Antonin Chambolle. Diagonal preconditioning for first order primal-dual algorithms in convex optimization. In *2011 IEEE International Conference on Computer Vision (ICCV)*, pages 1762–1769. IEEE, 2011.

[33] R Tyrrell Rockafellar. Convex analysis. princeton landmarks in mathematics, 1997.

[34] Leonid I Rudin, Stanley Osher, and Emad Fatemi. Nonlinear total variation based noise removal algorithms. *Physica D: Nonlinear Phenomena*, 60(1):259–268, 1992.

[35] Ivan W Selesnick, Ankit Parekh, and Ilker Bayram. Convex 1-d total variation denoising with non-convex regularization. *IEEE Signal Processing Letters*, 22(2):141–144, 2015.

[36] Emil Y Sidky, Rick Chartrand, John M Boone, and Xiaochuan Pan. Constrained minimization for enhanced exploitation of gradient sparsity: Application to CT image reconstruction. *IEEE Journal of Translational Engineering in Health and Medicine*, 2:1800418, 2014.

[37] Robert Tibshirani, Michael Saunders, Saharon Rosset, Ji Zhu, and Keith Knight. Sparsity and smoothness via the fused lasso. *Journal of the Royal Statistical Society: Series B (Statistical Methodology)*, 67(1): 91–108, 2005.

[38] Tuomo Valkonen. A primal–dual hybrid gradient method for nonlinear operators with applications to mri. *Inverse Problems*, 30(5):055012, 2014.

[39] Fenghui Wang, Zongben Xu, and Hong-Kun Xu. Convergence of bregman alternating direction method with multipliers for nonconvex composite problems. *arXiv preprint arXiv:1410.8625*, 2014.

[40] Fenghui Wang, Wenfei Cao, and Zongben Xu. Convergence of multi-block bregman admm for nonconvex composite problems. *arXiv preprint arXiv:1505.03063*, 2015.

[41] Huahua Wang and Arindam Banerjee. Bregman alternating direction method of multipliers. In *Advances in Neural Information Processing Systems*, pages 2816–2824, 2014.

[42] Huahua Wang, Arindam Banerjee, and Zhi-Quan Luo. Parallel direction method of multipliers. In *Advances in Neural Information Processing Systems*, pages 181–189, 2014.

[43] Yu-Xiang Wang, James Sharpnack, Alex Smola, and Ryan J Tibshirani. Trend filtering on graphs. *arXiv preprint arXiv:1410.7690*, 2014.

[44] Cun-Hui Zhang. Nearly unbiased variable selection under minimax concave penalty. *The Annals of Statistics*, 38(2):894–942, 2010.